\documentclass[12pt]{amsart}
\usepackage{amsthm,amsmath,amssymb,amscd,graphics,enumerate, stmaryrd,xspace,verbatim, epic, eepic,url,fullpage}

\usepackage{enumitem}

\usepackage[usenames,dvipsnames]{color}

\usepackage[colorlinks=true]{hyperref}
\usepackage{cleveref}

\usepackage[active]{srcltx}
\usepackage[all]{xypic}
\SelectTips{cm}{}
\usepackage{quiver}

   \newtheorem{theorem}{Theorem}[section]
      \newtheorem*{theorem*}{Theorem}
   \newtheorem{proposition}[theorem]{Proposition}

   \newtheorem{lemma}[theorem]{Lemma}

   \newtheorem{corollary}[theorem]{Corollary}
   
   \newtheorem*{conjecture*}{Conjecture}

\theoremstyle{definition}
          \newtheorem*{exercise*}{Exercise}
   
   \newtheorem{example}[theorem]{Example}
 
   \newtheorem*{example*}{Example}
   
   \newtheorem{definition}[theorem]{Definition}

   \newtheorem*{definition*}{Definition}
   
   \newtheorem{remark}[theorem]{Remark}

    \newtheorem{construction}[theorem]{Construction}

\newcommand{\RR}{{\mathbb{R}}}

\newcommand{\CC}{{\mathbb{C}}}
\newcommand{\QQ}{{\mathbb{Q}}}

\newcommand{\PP}{{\mathbb{P}}}
\newcommand{\ZZ}{{\mathbb{Z}}}

\def\Bl{{\rm Bl}}

\newcommand{\cA}{{\mathcal A}}
\newcommand{\cB}{{\mathcal B}}
\newcommand{\cC}{{\mathcal C}}
\renewcommand{\cD}{{\mathcal D}}

\newcommand{\cO}{{\mathcal O}}
\newcommand{\cP}{{\mathcal P}}

\newcommand{\fP}{{\mathfrak P}}

\newcommand{\nc}{{\operatorname{nc}}}

\def\<{\langle}
\def\>{\rangle}

\newcommand{\Spec}{\operatorname{Spec}}
\newcommand{\Specan}{\operatorname{Specan}}

\newcommand{\Proj}{\operatorname{Proj}}
\newcommand{\Projan}{\operatorname{Projan}}

\newcommand{\colim}{\operatorname{colim}}

\newcommand{\Hom}{{\operatorname{Hom}}}

\newcommand{\Aut}{{\operatorname{Aut}}}

\newcommand{\codim}{\operatorname{codim}}

\def\:{{\colon}}
\def\.{{,\dots,}}
\def\dim{{\rm dim}}

\def\Ex{\mathrm{Ex}}

\newcommand{\double}{\genfrac..{0pt}1
{\raise -1pt\hbox{$\scriptstyle\longrightarrow$}}{\raise 3pt\hbox
{$\scriptstyle\longrightarrow$}}}

\newcommand{\setmin}{\smallsetminus}

\renewcommand{\setminus}{\smallsetminus}

\renewcommand\H{\operatorname{H}}

\def\int{{\rm int}}

\def\tototi{\mathbin{\mathop{\otimes}\limits^{\raise-1pt\hbox
{$\scriptscriptstyle {\rm L}$}}}}

\def\indlim{\mathop{\vrule width0pt height7pt depth
4pt\smash{\lim\limits_{\raise 1pt\hbox to 14.5pt
{\rightarrowfill}}}}}
\def\projlim{\mathop{\vrule width0pt height7pt depth
4pt\smash{\lim\limits_{\raise 1pt\hbox to 14.5pt
{\leftarrowfill}}}}}

\newcommand\displaceamount{3pt}

\newcommand{\doubledown}{\ar@<\displaceamount>[d]\ar@<-\displaceamount>[d]}

\newcommand{\doubleup}{\ar@<\displaceamount>[u]\ar@<-\displaceamount>[u]}

\newcommand{\doubleright}{\ar@<\displaceamount>[r]\ar@<-\displaceamount>[r]}

\newcommand{\res}{{\operatorname{res}}}

\def\Xbar{\overline{X}}

\def\Supp{\mathrm{Supp\ \!}}
\def\Diff{\mathrm{Diff}}

\def\Bim{\mathrm{Bim}}

\def\phibar{\overline{\phi}}
\def\psibar{\overline{\psi}}
\def\Ybar{\overline{Y}}
\def\Zbar{\overline{Z}}

\def\Bs{\mathrm{Bs}}

\def\Min{\mathrm{Min}}

\definecolor{AL}{RGB}{10,200,100}

\numberwithin{equation}{section}

\begin{document}

\title{Semi-log canonical minimal model program for projective morphisms between complex analytic spaces}

\author{Alberto Landi}
\address{Department of Mathematics, Box 1917, Brown University,
Providence, RI, 02912, U.S.A}
\email{alberto\_landi@brown.edu}


\maketitle

\begin{abstract}
    We prove the existence of certain steps of the minimal model program for projective morphisms between complex analytic spaces in the semi-log canonical case, generalizing results of Ambro and Koll\'ar. This is obtained by providing a sufficient criterion for the finite generation of the relative log canonical algebra. As a consequence, we prove demi-normal analogues of results of Hacon--Xu and Birkar. As a further application, we derive the existence of crepant semi-dlt modifications of analytic slc pairs whose irreducible components are normal in codimension 1, extending a result of Hashizume.

    The methods are new even in the algebraic case, providing streamlined proofs of the algebraic counterparts to the statements proved in this work, without relying on Koll\'ar's gluing techniques.
\end{abstract}

\tableofcontents

\section{Introduction}

\subsection{Semi-log canonical pairs}
Given a projective \emph{log canonical (lc) pair} $(X,\Delta_X)$ and a projective morphism $f\colon X\rightarrow T$, finite generation of the relative \emph{log canonical algebra} $R_X:=R(X/T,K_X+\Delta_X)$ is a fundamental problem in higher-dimensional algebraic geometry. The minimal model program (MMP) provides powerful methods for establishing finite generation in many important cases, in particular for smooth varieties~\cite{BCHM,HM2,Cascini-Lazic_NewOutlookI}.
When $R_X$ is finitely generated and $K_X+\Delta_X$ is $f$-big, the induced birational map $X\dashrightarrow\Proj R_X$ yields a \emph{log canonical model over $T$}. The target carries a natural structure of an lc pair for which the log canonical divisor is ample over $T$; when $T$ is a point, such a pair is called \emph{stable}. Stable pairs form the basic objects parametrized by the \emph{Koll\'ar--Shepherd-Barron--Alexeev (KSBA) moduli space}~\cite{KSB88,Alexeev_MgnW_surfaces,Kollar-moduli}, which plays a central role in the moduli theory of higher-dimensional varieties.

As is well known in the case of curves, obtaining a compact moduli space requires allowing degenerations to reducible or self-intersecting varieties. This leads naturally to the notion of \emph{semi-log canonical} (\emph{slc}) pairs, or \emph{locally-stable pairs} \cite{Kollar-moduli}, introduced in~\cite{KSB88} for surfaces. The underlying variety of an slc pair is \emph{demi-normal}, that is, it satisfies Serre's $S_2$ condition and has at worst ordinary double points in codimension 1. A systematic treatment of slc pairs is substantially more difficult, as MMP steps cannot always be executed~\cite[Examples 4, 5]{Ambro-Kollar-slc-MMP} and finite generation of the log canonical ring can fail even for normal crossing surfaces without boundary~\cite{Kollar_examples_SNC}. Nevertheless, substantial progress has been made in the study of slc pairs since their introduction~\cite{FujinoAbundanceSLC,Kollar-singularities,Gongyo_abundance,fujino_theorems_slc,FG_fin_Brepr,Fujita_semi-terminal_modifications,Hacon_Xu_finBrepr_slc_abundance,Ambro-Kollar-slc-MMP,Hashizume_sdlt_models}.

\subsection{The analytic setting and main results}

In recent years, considerable effort has been devoted to extending the relative minimal model program to projective morphisms between complex analytic spaces~\cite{Fujino_MMPanalytic,Fujino_ConeContractionAnalytic,Enokizono-Hashizume-analytic-MMP,Fujino_RelativeLogPluricanonical,Fujino_QuasiLogAnalytic,Fujino_VanishingAnalytic,Enokizono_Hashizume_MMP_stacks}. Under suitable hypotheses (see \cite[page 3]{Fujino_MMPanalytic}), many results extend from the algebraic to the analytic setting with only minor modifications, although some technical obstacles remain. One important tool that the theory is currently missing is Koll\'ar's powerful gluing machinery \cite[\S5, \S9]{Kollar-singularities}, for which a suitable analogue is not currently available in the complex analytic setting. As a consequence, many results concerning analytic semi-log canonical pairs cannot be directly obtained from the corresponding results for normal pairs via gluing arguments. In particular, the following theorem is new in the complex setting, and its proof is new even in the algebraic framework. It constitutes the main result of the paper, and will provide our principal tool for studying the minimal model program for slc pairs in the analytic setting.

We will only work with complex analytic spaces that are \emph{Hausdorff} and \emph{second-countable}. We refer to \S\ref{subsec: singularities pairs} for background on singularities of pairs, including slc and semi-dlt pairs, as well as on log centers. See \S\ref{sec: good stable models} for the definition and properties of good stable models, which are a well-behaved demi-normal analogue of log canonical models.

\begin{theorem}\label{thm: general version}
    Let $(X,\Delta_X)$ be a semi-dlt complex analytic pair, $T$ a semi-normal complex analytic space, and $f\colon X\rightarrow T$ a projective morphism. Let $T^0\subset T$ be a dense Zariski-open subspace,
    and set $X^0=T^0\times_{T}X$, $\Delta_{X^0}=\Delta_X|_{X^0}$. Assume that $X^0$ intersects
    all log centers of $(X,\Delta_X)$ lying on the conductor $D_X$; in particular, $X^0$ intersects all strata of $D_X$.
    Suppose that $(X^0,\Delta_{X^0})$ admits a good stable model
    \[
    \begin{tikzcd}
        \phi^{0,c}\colon (X^0,\Delta_{X^0})\arrow[r,dashed] & (X^{0,c},\Delta_{X^{0,c}})
    \end{tikzcd}
    \]
    whose exceptional locus does not contain any stratum of $D_{X^0}$, and that the normalization pair $(\Xbar,D_{\Xbar}+\Delta_{\Xbar})$ admits a log canonical model over $T$.
    
    Then, $R_X:=R(X/T,K_X+\Delta_X)$ is locally finitely generated, and $\Projan_T R_X$ is the good stable model of $(X,\Delta_X)$ over $T$.
\end{theorem}

Note that the problem of finite generation of the log canonical ring of an slc pair is considerably more subtle than in the log canonical case: it can fail even for normal crossing algebraic surfaces without boundary \cite{Kollar_examples_SNC}, even though the log canonical ring of their normalization is known to be finitely generated. To address this problem, we assume the existence of a good stable model over a dense Zariski-open subset $T^0$ of $T$. This assumption is often satisfied both for families varying in moduli and when running the minimal model program, making it well suited to applications.

When the reduced part of the boundary is trivial, Theorem~\ref{thm: general version} admits an equivariant strengthening with respect to the group of B-bimeromorphic self-maps $\Bim(\overline{X}/T,D_{\Xbar}+\Delta_{\Xbar})$ on the normalization pair, as defined in \S\ref{subsec: B-bimeromorphic}.

\begin{theorem}\label{thm: equivariant sklt}
    With the notation of Theorem~\ref{thm: general version}, assume further that $\lfloor\Delta_X\rfloor=0$.
    Then, for every finite subgroup $G\subset\Bim(\Xbar/T,D_{\Xbar}+\Delta_{\Xbar})$, the inclusion
    \[
    \begin{tikzcd}
        R(X/T,K_X+\Delta_X)\cap R(\Xbar/T,K_{\Xbar}+D_{\Xbar}+\Delta_{\Xbar})^G\arrow[r,hookrightarrow] & R(\Xbar/T,K_{\Xbar}+D_{\Xbar}+\Delta_{\Xbar})
    \end{tikzcd}
    \]
    is a finite extension. In particular, $R_X:=R(X/T,K_X+\Delta_X)$ is locally finitely generated, and $\Projan_T R_X$ is the good stable model of $(X,\Delta_X)$ over $T$.
\end{theorem}

Theorem~\ref{thm: general version} allows us to deduce many results on finite generation and existence of MMP steps for slc pairs from the normal case. In particular, we are able to prove the following semi-log canonical version of~\cite[Theorem 1.1]{Hacon_Xu_dslt}.

\begin{corollary}\label{cor: Hacon-Xu generalization}
    Let $(X,\Delta_X)$ be a semi-dlt pair, $T$ a semi-normal scheme that is separated and essentially of finite type over $\CC$, and $f\colon X\rightarrow T$ a projective morphism. Let $T^0\subset T$ be a Zariski-open subspace,
    and set $X^0=T^0\times_{T}X$, $\Delta_{X^0}=\Delta_X|_{X^0}$. Assume that $X^0$ intersects all log canonical centers of $(X,\Delta_X)$, and all log centers lying on $D_X$; in particular, $X^0$ intersects all strata of the conductor $D_X$.
    Suppose that $(X^0,\Delta_{X^0})$ admits a good stable model
    \[
    \begin{tikzcd}
        \phi^{0,c}\colon (X^0,\Delta_{X^0})\arrow[r,dashed] & (X^{0,c},\Delta_{X^{0,c}})
    \end{tikzcd}
    \]
    over $T^0$ whose exceptional locus does not contain any stratum of $D_{X_0}$.
    
    Then, $R_X:=R(X/T,K_X+\Delta_X)$ is locally finitely generated, and $\Projan_T R_X$ is the good stable model of $(X,\Delta_X)$ over $T$.
\end{corollary}

The above result can be used to extend slc stable families from a punctured disk to the whole disk, by combining it with a generalization of semi-stable reduction to be developed in forthcoming joint work with Roberta Pagliaro, Simon Stojkovic, and Hao (Nick) Sun.

\subsection{Existence of MMP steps for analytic slc pairs}

We are interested in generalizing the existence of MMP steps for projective morphisms between normal complex analytic spaces to the non-normal case. In general, it is well-known that the existence of MMP steps for slc pairs can fail even in the algebraic case, see~\cite[Examples 4, 5]{Ambro-Kollar-slc-MMP} and~\cite[Example 5.4]{fujino_theorems_slc}. However, results of Fujino guarantee that the cone and contraction theorems hold for slc pairs, even in the complex analytic setting~\cite[Theorems 9.1, 9.2]{Fujino_QuasiLogAnalytic}.
We show that if the exceptional locus of a negative extremal contraction contains no log center lying on the conductor, we can then perform an MMP step. In the algebraic world, this has been proven in~\cite{Ambro-Kollar-slc-MMP,Hashizume_sdlt_models}. The first step is proving the following semi-log canonical version of~\cite[Theorem 1.2]{Enokizono-Hashizume-analytic-MMP}.

\begin{corollary}\label{cor: analog Birkar}
    Let $f\colon X\rightarrow T$ be a projective generically-finite morphism between demi-normal complex analytic spaces. Let $\Delta_X$ and $\Delta_X'$ be effective $\QQ$-divisors on $X$ such that $(X,\Delta_X+\Delta_X')$ is slc and $(X,\Delta_X)$ is semi-dlt. Suppose that $K_X+\Delta_X+\Delta_X'\sim_{\QQ,T}0$, and that $f$ is finite at very general points of every log center lying on $D_X$. Then, $R_X:=R(X/T,K_X+\Delta_X)$ is locally finitely generated, and $\Projan_T R_X$ is the good stable model of $(X,\Delta_X)$ over $T$.
\end{corollary}

As mentioned above, Corollary~\ref{cor: analog Birkar} allows us to construct certain MMP steps. Here, we have to be careful in asserting the existence of steps only over some open neighborhood of a Stein compact subset with a Noetherianity assumption, as in the main results of~\cite{Fujino_MMPanalytic,Enokizono-Hashizume-analytic-MMP}; see also Definition~\ref{def:slc-MMP steps}. The assumption on the compact set is often satisfied, see Theorem~\ref{thm: Siu} and Example~\ref{exm: property P4}.

\begin{corollary}\label{cor: existence slc MMP}
    Let $(X,\Delta_X)$ be an slc pair, $f\colon X\rightarrow T$ a projective morphism to a semi-normal Stein space $T$, and $W\subset T$ a Stein compact subset of $T$ such that $\Gamma(W,\cO_T)$ is Noetherian. Let $F\subset\overline{NE}(X/T;W)$ be a $(K_X+\Delta_X)$-negative extremal ray, and $\varphi=\varphi_F\colon X\rightarrow Z$ the corresponding contraction, defined after suitably restricting $T$ around $W$. Suppose that $\varphi$ is bimeromorphic and $\Ex(\varphi)\cap D_X$ contains no log centers. Then, $Z$ is demi-normal and there exists an MMP step around $W$
    \[
    \begin{tikzcd}
        (X,\Delta_X)\arrow[rr,dashed,"\psi"]\arrow[rd,"\varphi"'] & & (Y,\Delta_Y)\arrow[ld,"\varphi'"]\\
        & Z
    \end{tikzcd}
    \]
    Moreover, if $(X,\Delta_X)$ is semi-dlt, then so is $(Y,\Delta_Y)$.
\end{corollary}

To prove Corollary~\ref{cor: existence slc MMP}, we first assume $(X,\Delta_X)$ to be semi-dlt, in which case $(Y,\Delta_Y)$ is also semi-dlt, see Corollary~\ref{cor: existence semidlt MMP}. Then, we use the MMP for semi-dlt pairs to construct crepant semi-dlt modifications of analytic slc pairs whose irreducible components are normal in codimension 1, generalizing~\cite[Theorem 1.2]{Hashizume_sdlt_models} to the analytic case; see Theorem~\ref{thm: semidlt modifications} below. Finally, we use the semi-dlt modification (and the double cover of Remark~\ref{rmk: double cover}) to prove the general case of Corollary~\ref{cor: existence slc MMP}.

\begin{theorem}\label{thm: semidlt modifications}
    Let $(X,\Delta_X)$ be an slc pair such that the irreducible components of $X$ are normal in codimension 1, and $f\colon X\rightarrow T$ a projective morphism to a complex analytic space. Let $W\subset T$ be a Stein compact subset of $T$ such that $\Gamma(W,\cO_T)$ is Noetherian. Then, after suitably restricting $T$ around $W$, there exists a semi-dlt pair $(Y,\Delta_Y)$ projective over $T$, and a projective bimeromorphic morphism $g\colon Y\rightarrow X$ over $T$, such that
    \begin{enumerate}
        \item $g$ is an isomorphism over a Zariski-open subset of $X$ intersecting all components of the conductor of $D_X$;
        \item $g^*(K_X+\Delta_X)=K_Y+\Delta_Y$, that is, $g$ is crepant;
        \item $K_Y$ is $\QQ$-Cartier;
        \item $(Y,0)$ is semi-terminal in the sense of Fujita~\cite{Fujita_semi-terminal_modifications}.
    \end{enumerate}
\end{theorem}

We remark that our notions of dlt and semi-dlt pair differ slightly from those introduced in~\cite{Fujino_MMPanalytic,Enokizono-Hashizume-analytic-MMP}; see Definition~\ref{def: dlt/sdlt} and Remark~\ref{rmk: discussion on dlt and semidlt}.

Although requiring control on the log centers lying on the conductor seems a strong condition, \cite[Examples 4, 5]{Ambro-Kollar-slc-MMP} and~\cite[Example 5.4]{fujino_theorems_slc} show that these conditions are needed and the results of this paper are at least close to optimal.

\subsection{The algebraic case}\label{subsec: discussion algebraic case}

The proofs of all the results above carry over, with only minor modifications, to the setting of separated schemes essentially of finite type over an algebraically closed field of characteristic zero. Moreover, in the algebraic setting, the results hold globally over the base $T$, and many subtleties related to the lack of Noetherianity of complex analytic spaces are absent. In particular, this work provides alternative proofs to the results of Ambro and Koll\'ar in \cite{Ambro-Kollar-slc-MMP}, and Hashizume in \cite{Hashizume_sdlt_models}; these proofs entirely avoid the use of Koll\'ar's gluing theory. While powerful, the gluing machinery is rather complicated and, so far, seems less amenable to generalization to the complex analytic setting.
The proofs in the language of schemes are available in the preprint~\cite{Landi_FiniteGenerationSemidlt}, which is not intended for publication and has been superseded by the present paper.

The algebraic versions of the results in this paper will be used in forthcoming work with Roberta Pagliaro, Simon Stojkovic, and Hao (Nick) Sun to give a streamlined proof of the valutative criterion of properness for the KSBA of stable pairs \cite{KSB88,Alexeev_MgnW_surfaces,Hacon_Xu_dslt,Kollar-moduli}.

\subsection{Structure of the paper}

In \S\ref{sec: preliminaries}, we recall the main definitions and results about complex analytic spaces, their singularities, and their minimal model program. Good stable models for slc pairs are introduced and studied in \S\ref{sec: good stable models}. In \S\ref{sec: finite generation}, we prove Theorems~\ref{thm: equivariant sklt} and~\ref{thm: general version} in this order. Moreover, in \S\ref{subsec: corollaries} we prove Corollaries~\ref{cor: Hacon-Xu generalization} and~\ref{cor: analog Birkar}. We leverage these results in \S\ref{sec: slc MMP} to establish the existence of certain MMP steps for projective morphisms of complex analytic pairs with slc singularities. First, we treat the semi-dlt case in \S\ref{subsec: sdlt MMP}, which is used in \S\ref{subsec: sdlt modification} to establish the existence of crepant semi-dlt modifications of slc pairs whose irreducible components are normal in codimension 1. Using this, we can extend the existence of certain MMP steps to slc pairs in \S\ref{subsec: slc MMP}.

\subsection{Conventions and Notation}

Throughout the paper, every complex analytic space is assumed to be Hausdorff and second-countable.

Let $f\colon X\rightarrow T$ be a projective morphism between complex analytic spaces with $X$ demi-normal, and let $D$ be a $\QQ$-Cartier $\QQ$-divisor on $X$. For every $m\geq0$, we denote by $\cO_X(mD)$ the reflexive sheaf $\cO_X(\lfloor mD\rfloor)$, and set
\[
    R(X/T,D):=\oplus_{m\geq0}f_*\cO_X(mD),
\]
which has a natural structure of $\cO_T$-algebra. When $D$ and $T$ are clear from context, we simply write $R_X$ in place of $R(X/T,D)$.

Given an integer $\ell>0$, a complex analytic space $T$, and a sheaf of graded $\cO_T$-algebras $\cA=\oplus_{n\geq0} \cA_n$, we denote by $\cA^{(\ell)}:=\oplus_{n\geq0}\cA_{\ell n}$ the $\ell$-th \emph{Veronese subalgebra} of $\cA$.

\subsection*{Acknowledgments}

This project was initially motivated by a question raised by Dan Abramovich about simplifying the proof of the properness of the KSBA moduli space of stable pairs. I am grateful to him for his suggestions and support. I am also indebted to Christopher Hacon, who suggested generalizing the methods to the complex analytic case, and patiently answered several questions. I would also like to thank Alessio Corti, Osamu Fujino, Kenta Hashizume, and J\'anos Koll\'ar for generously answering questions and providing feedback on previous versions of the manuscript. Finally, I thank Brendan Hassett, Eric Jovinelly, Roberta Pagliaro, Simon Stojkovic, and Hao (Nick) Sun for useful conversations on this and related projects.
This research is supported by NSF grant DMS-2401358.

\subsection*{AI disclosure}

To prepare the previous version of the manuscript \cite{Landi_FiniteGenerationSemidlt}, concerning only the algebraic case, the author used free versions of ChatGPT and Gemini in a limited capacity for proofreading, improving the clarity and exposition of this paper, and occasional review of proofs. This review did not result in modifications to the proofs. The mathematical results, arguments, and ideas were developed independently by the author. The writing of the present version of the paper was completely performed by the author with no use of AI.

\section{Preliminaries}\label{sec: preliminaries}

In this section we will recall basic notions and properties that we will need in the rest of the paper, and prove some preliminary results.

\subsection{Complex Analytic Spaces}

Compared to schemes of finite type over $\CC$, complex analytic spaces hide many difficulties, most of which are related to the lack of Noetherianity of rings of global sections. In this subsection we will highlight the main differences and cover the results needed for the rest of the paper. For the basic theory of complex analytic spaces, we refer to \cite{Fischer_BookAnalytic,Nakayama_ZariskiDecomposition,Fujino_MMPanalytic,Enokizono-Hashizume-analytic-MMP,Das_Hacon_Kahler3foldsMMP,Das_Hacon_Pau_4dimensionKahler}.

Every complex analytic space is assumed to be \emph{Hausdorff} and \emph{second-countable}. A \emph{complex manifold} is a smooth complex analytic space.

\subsubsection{Stein spaces and Stein compact subsets}\label{subsec: Stein spaces}

Stein spaces are the complex analytic counterpart to affine schemes. Recall that a complex analytic space $X$ is \emph{Stein} if the global sections of $\cO_X$ separate points, generate the maximal ideal of $\cO_{X,x}$ for every $x\in X$, and $X$ is \emph{holomorphically convex}; see~\cite[Definitions 2.1, 2.2]{Fujino_MMPanalytic} for this last notion.
For instance, the \emph{analytification} of an affine scheme is Stein.
Moreover, any complex analytic space is locally a Stein space.

The most important properties of Stein spaces are that any coherent sheaf is globally generated (Cartan's Theorem A \cite[Theorem 0.37 A]{Fischer_BookAnalytic}) and has no higher cohomology (Cartan's Theorem B \cite[Theorem 0.37 B]{Fischer_BookAnalytic}). As for affine schemes, the vanishing of the first cohomology group for any coherent sheaf actually characterizes Stein spaces, see~\cite[Lemma 2.7]{Fujino_MMPanalytic}.

Even though we saw that Stein spaces satisfy many properties that resemble affine schemes, they do not enjoy nice finiteness properties. Specifically, the ring of global sections $H^0(X,\cO_X)$ is not Noetherian for any positive dimensional Stein space, see~\cite[Example 2.9]{Fujino_MMPanalytic}.
To solve this issue we have to look at holomorphic functions defined around a fixed compact subset. A particularly convenient class of compact subsets consists of \emph{Stein compact} subsets.

\begin{definition}\label{def: Stein compact}
    A compact subset $K$ of a complex analytic space is called \emph{Stein compact} if it admits a fundamental system of Stein open neighborhoods.
\end{definition}

Note that Stein compact subsets are not compact subsets that are Stein, as these would automatically be 0-dimensional by the maximum principle.
Moreover, it is possible to construct many Stein compact subsets by taking holomorphic convex hulls of compact subset of Stein spaces, see~\cite[Lemma 2.5]{Fujino_MMPanalytic}.

Being a Stein compact subset is not enough to guarantee nice finiteness properties.

\begin{definition}\label{def: P4 property}
    Let $X$ be a complex analytic space, and $K$ be a compact subset of $X$. We say that $K$ satisfies \emph{property $\cP$}, if the intersection $K\cap Z$ has only finitely many connected components for every analytic subset $Z$ defined on an open neighborhood of $K$.
\end{definition}

The relevance of the above definition is explained by the following theorem of Siu.

\begin{theorem}[{\cite[Theorem 1]{Siu_Noetherianity}}]\label{thm: Siu}
    Let $K$ be a Stein compact subset of a complex analytic space $X$. Then, $K$ satisfies property $\cP$
    if and only if
    \[
        \cO_X(K):=\Gamma(K,\cO_Y)=\colim_{K\subset U}\Gamma(U,\cO_X)
    \]
    is Noetherian, where $U$ runs through all the open neighborhoods of $K$.
\end{theorem}

\begin{example}\label{exm: property P4}
    For instance, for every point $p\in X$ of a complex analytic space the stalk $\cO_{X,x}$ is Noetherian. Moreover, any Stein compact subset that is semi-analytic satisfies property $\cP$ \cite[Corollary 2.7]{Bierstone_Milman_semianalytic}. As a non-example, the Cantor set is a Stein compact subset of $\CC$ that does not satisfy property $\cP$, see~\cite[Example 2.12]{Fujino_MMPanalytic}.
\end{example}

\subsubsection{Finite generation of sheaves of algebras}

Following Fujino~\cite{Fujino_MMPanalytic}, we recall the basic notion of finite generation for a sheaf of algebras, and some properties.

\begin{definition}[{\cite[Definition 2.23]{Fujino_MMPanalytic}}]\label{def: finitely generated algebra}
    Let $X$ be a complex analytic space. A graded $\cO_X$-algebra $\cA=\oplus_{m\geq0}\cA_m$ is called \emph{finitely generated} if there exists a surjective morphism
    \[
    \begin{tikzcd}
        \cO_X[x]:=\cO_X[x_1,\ldots,x_{\ell}]=\cO_X\otimes\CC[x]\arrow[r,"\varphi",twoheadrightarrow] & \cA
    \end{tikzcd}
    \]
    for some $\ell$, where $x_i$ is mapped to an homogeneous element. If further the kernel of $\varphi$ is generated by a finite number of polynomials in $\H^0(X,\cO_X)[x]$, then we say that $\cA$ is \emph{of finite presentation}. We say that $\cA$ is \emph{locally finitely generated} (resp. \emph{locally of finite presentation}) if there exists an open cover $X=\cup_{\lambda}U_{\lambda}$ such that $\cA|_{U_\lambda}$ is finitely generated (resp. of finite presentation).
\end{definition}

We will always work with $\cO_X$-algebras $\cA=\oplus_{m\geq0}\cA_m$ where $\cA_m$ are coherent $\cO_X$-modules. In this case, the notion of finite generation and finite presentation agree.

\begin{lemma}[{\cite[Lemma 2.24]{Fujino_MMPanalytic}}]\label{lem: finite presentation}
    Let $X$ be a complex analytic space and let $\cA=\oplus_{m\geq0}\cA_m$ be a locally finitely generated graded $\cO_X$-algebra such that $\cA_m$ are all coherent $\cO_X$-modules. Then, $\cA$ is locally of finite presentation.
\end{lemma}

In the proofs of Theorems~\ref{thm: general version} and~\ref{thm: equivariant sklt}, we will often make use of the following lemma, that reduces the finite generation of a $\cO_X$-algebra $\cA$ to the finite generation of its $\ell$-Veronese subalgebra $\cA^{(\ell)}$.

\begin{lemma}\label{lem: veronese finite generation}
    Let $X$ be a complex analytic space and let $\cA=\oplus_{m\geq0}\cA_m$ be a graded $\cO_X$-algebra such that $\cA_m$ is a coherent $\cO_X$-module for every $m$ and $\cA(U)$ is a reduced algebra for all open subsets $U\subset X$.
    Let $\ell$ be any positive integer. Then, the graded $\cO_X$-algebra $\cA$ is locally finitely generated if and only if so is $\cA^{(\ell)}$.
\end{lemma}
\begin{proof}
    If $\cA(U)$ is an integral domain for every connected relatively compact open $U\subset X$, the result follows from~\cite[Lemma 2.26]{Fujino_MMPanalytic}. Moreover, the same argument as in the proof of loc.\! cit. shows that, in the general case, if $\cA$ is locally finitely generated then so is $\cA^{(\ell)}$.
    
    From now on, let us assume that $\cA^{(\ell)}$ is locally finitely generated, and fix $p\in X$. We can choose a connected relatively-compact Stein open neighborhood $U\subset X$ of $p$, and a connected Stein compact subset $K\subset X$ satisfying $\cP$ such that $p\in K\subset U\subset X$. Recall that the property $\cP$ is equivalent to $\Gamma(K,\cO_X)$ being Noetherian by Siu's Theorem~\ref{thm: Siu}, while the Stein property of $U$ implies that $\cA_{m,x}$ is generated by $\cA_m(U)$ for all $m\geq0$ and $x\in U$. Eventually after restricting $K$ and $U$, we can further assume that $\cA^{(\ell)}$ is finitely generated over $U$.
    Choose a surjection
    \[
    \begin{tikzcd}
        \cO_U[x]:=\cO_U[x_1,\ldots,x_r]\arrow[r,twoheadrightarrow] & \cA|_U^{(\ell)}
    \end{tikzcd}
    \]
    and note that $\cO_U(K)[x]\rightarrow\cA^{(\ell)}(K)$ is still surjective by~\cite[Lemma 2.19]{Fujino_MMPanalytic}, since both $\cO_U[x]$ and $\cA|_U^{(\ell)}$ admit $\cO_U$-exhaustions (\cite[Definition 2.18]{Fujino_MMPanalytic}).
    Since $\cO_X(K)$ is Noetherian, $\cA^{\ell}(K)$ is itself Noetherian hence its set of minimal primes is finite. Since minimal primes of a graded ring are homogeneous by~\cite[\href{https://stacks.math.columbia.edu/tag/00JU}{Tag 00JU}]{stacks} and the set of homogeneous primes of $\cA(K)$ and $\cA^{(\ell)}(K)$ coincide, the set of minimal primes $\Min(\cA)$ of $\cA(K)$ coincides with that of $\cA^{(\ell)}(K)$, hence it is finite. For every $\fP\in\Min(\cA)$, we set $\fP^{(\ell)}=\fP\cap\cA^{(\ell)}(K)$.
    Now, $\cA^{(\ell)}(K)/\fP^{(\ell)}$ is a finitely generated graded $\cO_X(K)$-algebra and coincides with $(\cA(K)/\fP)^{(\ell)}$. Since $\cO_X(K)$ is Noetherian and $\cA^{(\ell)}(K)/\fP^{(\ell)}$ is an integral domain, it follows that $\cA(K)/\fP$ is also finitely generated over $\cO_X(K)$, hence $\cA^{(\ell)}(K)/\fP^{(\ell)}\subset\cA(K)/\fP$ is a finite extension; see the proof of~\cite[Lemma 2.26]{Fujino_MMPanalytic}. Since $\cA(K)$ injects into the finite $\cA^{(\ell)}(K)$-module $\prod_{\fP\in\Min(\cA)}\cA(K)/\fP$ and $\cA^{(\ell)}(K)$ is Noetherian, it follows that $\cA(K)$ is also a finite $\cA^{(\ell)}(K)$-module. In particular, $\cA(K)$ is also a finitely generated $\cO_X(K)$-algebra. Since $\cA_{m,x}$ is generated by $\cA_m(K)$ for all $m\geq0$ and $x\in K$, we have that $\cA$ is finitely generated around $p$. This concludes.
\end{proof}

We will apply the above lemma to the algebra $\cA:=R(X/T,D)$, where $X$ is a reduced complex analytic space, $D$ a $\QQ$-Cartier divisor on $X$, and $f\colon X\rightarrow T$ a projective morphism to a complex analytic space.

The same argument as in Lemma~\ref{lem: veronese finite generation} proves the following version of the Artin-Tate lemma.

\begin{lemma}[Artin-Tate Lemma]\label{lem: Artin-Tate}
    Let $X$ be a complex analytic space, and $\cA\hookrightarrow\cB$ an injective morphism of graded $\cO_X$-algebras whose summands $\cA_m$ and $\cB_m$ are coherent $\cO_X$-modules for all $m\geq0$. Suppose that $\cB$ is locally finitely generated, and that $\cB$ is locally a finite $\cA$-module. Then, $\cA$ is locally finitely generated.
\end{lemma}
\begin{proof}
    By working locally on $X$, we can assume $X$ to be Stein and the algebras $\cA$ and $\cB$ to be globally generated by Cartan's Theorem A. We can also assume that there exists a surjection $\cA^{\oplus \ell}\twoheadrightarrow\cB$ of $\cA$-algebras. Given $x\in X$, let $K$ be a Stein compact neighborhood of $x$ satisfying property $\cP$. Then, $\cA(K)^{\oplus\ell}\rightarrow\cB(K)$ is again surjective by~\cite[Lemma 2.19]{Fujino_MMPanalytic}. Since $\cO_X(K)$ is Noetherian, the classical Artin-Tate Lemma~\cite[\href{https://stacks.math.columbia.edu/tag/00IS}{Tag 00IS}]{stacks} implies that $\cA(K)$ is finitely generated over $\cO_X(K)$. Since $\cA(K)$ generates $\cA$ at all points of $K$, the results follows.
\end{proof}

We will also need the following simple lemma.

\begin{lemma}\label{lem: pullback integral extension}
    Let $X$ be a complex analytic space, and $\cA$, $\cB$, $\cC$ be graded $\cO_X$-algebras whose summands $\cA_m$, $\cB_m$, $\cC_m$ are coherent $\cO_X$-modules for all $m\geq0$. Let $\cA\hookrightarrow\cB$ be an integral injective homomorphism, and $\cC\twoheadrightarrow\cB$ a surjective homomorphism. Then, $\cD:=\cA\times_{\cB}\cC\hookrightarrow\cC$ is again integral.
\end{lemma}
\begin{proof}
    This follows immediately from the analogous statement between rings: if $A\hookrightarrow B$ is an integral extension of rings and $C\twoheadrightarrow B$ a surjective homomorphism of rings, then $D:=A\times_BC\hookrightarrow C$ is an integral extension. To prove this last statement, let $c\in C$ and $b$ its image in $B$. By assumption, there exists a monic polynomial $p(x)\in A[x]$ such that $p(b)=0$. Since $D\rightarrow A$ is again surjective, there exists a monic polynomial $q(x)\in D[x]$ whose image along $D[x]\rightarrow A[x]$ is $p(x)$. It follows that $q(c)\in\ker(C\rightarrow B)\subset D$, hence $Q(x):=q(x)-q(c)\in D[x]$ is the desired monic polynomial such that $Q(c)=0$.
\end{proof}

\subsubsection{$\Specan$ and $\Projan$}

Given a complex analytic space $X$ and a $\cO_X$-algebra $\cA$ locally of finite presentation, we can construct the \emph{relative analytic spectrum} $\Specan_X\cA\rightarrow X$ satisfying the usual universal property; see~\cite[II.1.4]{Nakayama_ZariskiDecomposition}, or the more general construction in~\cite{Hakim_book_Spec}.
As in the algebraic case, the analytic spectrum allows us to construct the Stein factorization of proper morphisms, (semi-)normalizations, and to establish Zariski's Main Theorem \cite{Fischer_BookAnalytic,Nakayama_ZariskiDecomposition}.

Similarly, given a graded $\cO_X$-algebra $\cA=\oplus_{m\geq0}\cA_m$ locally of finite presentation, we can construct the \emph{relative analytic Proj} $\Projan_X\cA\rightarrow X$, satisfying the usual universal property; see~\cite[Lemma 1.8]{Nakayama_ZariskiDecomposition}.

We refer to~\cite[\S1]{Nakayama_ZariskiDecomposition} for more details and the rigorous statements that show how $\Specan$ and $\Projan$ have the same properties as their algebraic counterparts.

\subsubsection{Meromorphic maps}\label{subsec: meromorphic maps}

The complex analytic counterpart to rational (resp. birational) maps are \emph{meromorphic} (resp. \emph{bimeromorphic}) maps. We recall their definition for the reader convenience. We refer to~\cite[\S II.1.d]{Nakayama_ZariskiDecomposition} for details.

\begin{definition}\label{def: meromorphic maps}
    A \emph{meromorphic map} $f\colon X\dashrightarrow Y$ of complex analytic spaces is defined by the graph $\Gamma_f\subset X\times Y$ such that $\Gamma_f$ is a subspace of $X\times Y$ and that the first projection is proper and an isomorphism over a Zariski-open dense subset of $X$. If further the second projection $\Gamma_f\rightarrow Y$ is proper, we say that $f$ is proper. If $f$ is a proper meromorphic map and $\Gamma_f\rightarrow Y$ is an isomorphism over a Zariski-open dense subset of $Y$, then $f\colon X\dashrightarrow Y$ is called a \emph{bimeromorphic map}. If $f\colon X\dashrightarrow Y$ is a proper bimeromorphic map, the \emph{exceptional locus} $\Ex(f)$ of $f$ is the complement of the largest open subset $U\subset X$ such that $f$ is a morphism over $U$ and $f|_U:U\rightarrow f(U)$ is an isomorphism. A \emph{bimeromorphic contraction} is a bimeromorphic map $f\colon X\dashrightarrow Y$ such that $\Ex(f^{-1})$ has codimension at least 2. If both $f$ and $f^{-1}$ are bimeromorphic contractions, we say that $f$ is \emph{small}.
\end{definition}

Note that our definition of exceptional locus differs slightly from~\cite[Definition 2.30]{Fujino_MMPanalytic}, and coincides with the complement of the largest Zariski-open $U\subset X$ over which $\pi_1:\Gamma_f\rightarrow X$ is an isomorphism and $\pi_1^{-1}(U)\rightarrow Y$ is an isomorphism onto its image. Moreover, $\Ex(f)$ is analytic.

\subsection{Singularities of Complex Analytic Spaces}

In this subsection we define the main classes of singularities that we will consider and briefly recall some of their properties. The definitions of dlt and semi-dlt pairs depart slightly from the recent literature on singularities of complex analytic spaces, thus we will explain the differences.

We refer to~\cite[\S2, \S3]{Fujino_MMPanalytic} for the definition of $\QQ$-divisors on complex analytic varieties, their equivalence relations, and the various operations that one can perform on them, like the round-down. The main difference is that a $\QQ$-divisor is only assumed to be a \emph{locally finite} sum of prime divisors with rational coefficients, instead of a finite sum. Note that the sum is finite over any relatively compact open subset, so this will not present an issue.

The notions of \emph{bigness}, \emph{nefness}, \emph{ampleness}, and \emph{semi-ampleness} are similar to the algebraic case, and coincide with the usual notions when working over relatively compact open subsets. Since the small differences of these concepts between the algebraic and complex case will not play a role in this paper, we limit ourselves to refer to~\cite[\S2]{Fujino_MMPanalytic} for details.

\subsubsection{Canonical divisor class and discrepancies}\label{subsec: canonical class}

Before introducing the main types of singularities we will work with, we briefly discuss the notion of canonical divisor in the complex analytic setting; we refer to~\cite[\S4]{Nakayama_ZariskiDecomposition} and~\cite{Fujino_MMPanalytic,Enokizono-Hashizume-analytic-MMP,Fujino_QuasiLogAnalytic} for the details.

First, recall that a reduced complex analytic space $X$ is \emph{demi-normal} if it is nodal in codimension 1 and satisfies Serre's property $S_2$.
Note that the open locus $X^{\nc}$ of $X$ consisting of smooth or normal crossing points carries an invertible canonical sheaf $\omega_{X^{\nc}}$. The \emph{canonical sheaf} $\omega_X$ of $X$ is the unique reflexive sheaf whose restriction to $X^{\nc}$ coincides with $\omega_{X^{\nc}}$.

Contrarily to the algebraic case, in general there is no global Weil divisor $K_X$ such that $\cO_X(K_X)\simeq\omega_X$. However, such divisor always exists locally on $X$. In particular, given any $\QQ$-divisor $\Delta_X$ on $X$ and $x\in X$, we say that $K_X+\Delta_X$ is $\QQ$-Cartier at $x$ if there exists an open neighborhood $U_x$ of $x$ and a Weil divisor $K_{U_x}$ on $U_x$ that represents $\omega_{U_x}$, such that the $\QQ$-divisor $K_{U_x}+\Delta_{X}|_{U_x}$ is $\QQ$-Cartier. We say that $K_X+\Delta_X$ is $\QQ$-cartier if it is $\QQ$-Cartier at all points.

We will use the symbol $K_X$ as a formal divisor class with
an isomorphism $\cO_X(K_X)\simeq\omega_X$ and call it the \emph{canonical divisor} of $X$, if there is no danger of confusion.

\begin{remark}\label{rmk: existence KX}
In many important situations, $K_X$ exists as an actual Weil divisor.

Given a proper bimeromorphic map $f\colon X\rightarrow Y$ between demi-normal complex analytic spaces with $f_*\cO_X\simeq\cO_Y$, for any small Stein open $U\subset Y$ on which there is a well-defined canonical divisor $K_U$, it is possible to find a canonical divisor $K_{f^{-1}U}$ whose pushforward along $f|_U$ is $K_U$. This allows us to define \emph{discrepancies} and \emph{minimal log discrepancies} in the usual way; see~\cite{Nakayama_ZariskiDecomposition,Fujino_MMPanalytic,Enokizono-Hashizume-analytic-MMP} for more details.

Another situation where $K_X$ can be defined is when $f\colon X\rightarrow T$ is a projective morphism obtained by restriction to a relatively compact Stein open subset of the base (\cite[Remark 2.43]{Fujino_MMPanalytic}). To see this, let $L$ be a $\pi$-ample line bundle on $X$. Then, both $H^0(X,L^{\otimes m})$ and $H^0(X,\omega_X\otimes L^{\otimes m})$ are non-zero for $m$ sufficiently large by~\cite[Theorem 2.12]{Das_Hacon_Pau_4dimensionKahler}, and the existence of $K_X$ follows. Since the statements of Theorems~\ref{thm: general version} and~\ref{thm: equivariant sklt} are local on the base, in the case we are interested in we will be able to treat $K_X$ as an honest divisor; for this reason, we will ignore further subtleties.
\end{remark}

\subsubsection{Singularities of pairs}\label{subsec: singularities pairs}

Recall that a reduced and equi-dimensional complex analytic space is \emph{demi-normal} if it $S_2$ and is nodal in codimension 1.
Given a demi-normal complex analytic space $X$, we denote by $D_X$ its \emph{conductor}, which can be described as the closure of the nodal locus.

A \emph{normal pair} (resp. \emph{demi-normal pair}) $(X,\Delta_X)$ consists of a normal (resp. demi-normal) complex analytic space $X$ and an effective $\QQ$-divisor $\Delta_X$ (resp. not containing any component of the conductor $D_X$ of $X$) such that $K_X+\Delta_X$ is $\QQ$-Cartier. We will also assume that $\Delta_X$ is a \emph{boundary}, that is, the prime divisors appearing in the sum defining $\Delta_X$ have coefficient at most 1.

The definitions of \emph{log canonical} (\emph{lc}), \emph{Kawamata log terminal} (\emph{klt}), \emph{purely log terminal} (\emph{plt}) and \emph{semi-log canonical} (\emph{slc}) singularities are completely analogous to the algebraic case; see for instance~\cite[\S3]{Fujino_MMPanalytic} or~\cite[\S2.5]{Enokizono-Hashizume-analytic-MMP}. Similarly, for an lc pair $(X,\Delta_X)$ we can define \emph{log centers} as the centers of divisors over $X$ with negative discrepancy, and \emph{log canonical centers} (or \emph{lc centers}) as the centers of divisors over $X$ with discrepancy $-1$. For an slc pair, log centers (resp. log canonical centers) can be defined analogously, or equivalently as the images of log centers (resp. log canonical centers) of the normalization pair $(\Xbar,D_{\Xbar}+\Delta_{\Xbar})$ (\S\ref{subsec: normal and demi-normal modification}) along the normalization map $\pi:\Xbar\rightarrow X$.

We write down the definition for dlt and semi-dlt singularities, as the definition we adopt differs slightly from the literature on the minimal model program for complex analytic spaces. As usual, \emph{SNC} stands for \emph{simple normal crossing}.

\begin{definition}\label{def: SNC and semi-SNC pairs}
    A normal pair $(X,\Delta_X)$ is an \emph{SNC pair} if $X$ is smooth and $\Delta_X$ is SNC. A demi-normal pair $(X,\Delta_X)$ is \emph{semi-SNC} if $X$ is SNC and $\Delta_X+D_X$ is SNC, with $\Delta_X$ and $D_X$ having no common component. Equivalently, Zariski-locally on $X$ there exists an embedding $(X,\Delta_X)\hookrightarrow(Y,X+\Delta_Y)$ into an SNC pair with $\Delta_Y|_X=\Delta_X$.
\end{definition}

\begin{definition}\label{def: dlt/sdlt}
    A normal pair $(X,\Delta_X)$ is \emph{divisorial log terminal} (\emph{dlt}) if there exists an open subspace $U\subset X$ whose complement $Z=X\setminus U$ is analytic, $(U,\Delta_U)$ is SNC, and $a(E,X,\Delta_X)>-1$ for every irreducible divisor $E$ over $X$ with $\mathrm{center}_X(E)\subset Z$.

    A demi-normal pair $(X,\Delta_X)$ is \emph{semi-divisorial log terminal} (\emph{semi-dlt}) if there exists an open subspace $U\subset X$ whose complement $Z=X\setminus U$ is analytic, $(U,\Delta_U)$ is semi-SNC, and $a(E,X,\Delta_X)>-1$ for every irreducible divisor $E$ over $X$ with $\mathrm{center}_X(E)\subset Z$.
\end{definition}

\begin{remark}\label{rmk: analycity SNC locus}
    Let $(X,\Delta_X)$ be a demi-normal pair. Then, the complement of any of the following open subsets is analytic:
    \begin{enumerate}
        \item\label{point: singular locus is analytic} the locus where $X$ is smooth;
        \item\label{point: non-SNC locus is analytic} the locus where $X$ is SNC;
        \item the locus where $(X,\Delta_X)$ is SNC;
        \item the locus where $(X,\Delta_X)$ is semi-SNC;
    \end{enumerate}
    In particular, in Definition~\ref{def: dlt/sdlt} the fact that we require $X\setminus U$ to be analytic is not restrictive. Point~\eqref{point: singular locus is analytic} follows from~\cite[Proposition 10.24]{Remmert_ComplexAnalyticSpaces}. In proving~\eqref{point: non-SNC locus is analytic} we can assume $X$ to be pure dimensional, and by point~\eqref{point: singular locus is analytic} and~\cite[10.22]{Remmert_ComplexAnalyticSpaces}, we can assume that all irreducible components of $X$ are smooth and the embedding dimension at each point to be at most $\dim\ \! X+1$. Since the problem is local, we can also assume $X$ to be embedded as a divisor in a complex analytic manifold, and that each irreducible component $X_i$ is defined by one equation $f_i$. The requirement that the $f_i$'s define a regular sequence at a point is open analytic. The last two points follow from~\eqref{point: singular locus is analytic} and~\eqref{point: non-SNC locus is analytic}.
\end{remark}

\begin{remark}\label{rmk: discussion on dlt and semidlt}
    In~\cite{Fujino_MMPanalytic,Enokizono-Hashizume-analytic-MMP} and other texts, the notion of dlt is more restrictive, so we report it here with a different name. We say that a normal pair $(X,\Delta_X)$ is \emph{F-dlt} if there exists a proper bimeromorphic morphism $f\colon Y\rightarrow X$ from a smooth complex variety $Y$ such that $\Ex(f)\cup\Supp f_*^{-1}\Delta_X$ are simple normal crossing divisors on $Y$ and that the discrepancy coefficient $a(E,X,\Delta_X)>-1$ holds for every $f$-exceptional divisor $E$. In the literature, F-dlt pairs are simply called dlt; we add the letter F to distinguish them from Definition~\ref{def: dlt/sdlt} and in honor of Fujino. First of all, we note that any F-dlt pair is also dlt. Indeed, if $U\subset X$ is the maximal open over which $f$ is an isomorphism, the pair $(U,\Delta_X|_U)$ is SNC and the discrepancy of any divisor over $X$ whose center lies on the complement of $U$ is strictly greater than $-1$. Moreover, $X\setminus U$ is analytic (see also Remark~\ref{rmk: analycity SNC locus}). In general, a dlt pair can fail to be F-dlt because log resolutions are only guaranteed to exist over relatively compact subsets, see for instance~\cite[Theorem 2.16]{Das_Hacon_Kahler3foldsMMP}. More precisely, \cite[Lemma 3.9]{Fujino_MMPanalytic} shows that any dlt pair is F-dlt when restricted to any relatively compact open subset; note that \cite[Lemma 3.9]{Fujino_MMPanalytic} is stated for F-dlt pairs but the proof actually uses our weaker dlt notion. An advantage of our definition is that every klt pair is dlt, while it is F-dlt only when restricted to relatively compact open subsets.

    For semi-dlt pairs, the situation is slightly different. In \cite[Definition 2.6]{Fujino_RelativeLogPluricanonical}, semi-dlt is defined as slc pairs whose irreducible components are normal and F-dlt when endowed with the induced boundary. If we ignore the small difference between dlt and F-dlt pairs, our notion of semi-dlt pair is strictly more stringent, even in the algebraic case (\cite[Definition 5.19]{Kollar-singularities}).
\end{remark}

All klt pairs are particular examples of dlt pairs, and the converse holds if $\lfloor\Delta_X\rfloor=0$. This justifies the following definition.

\begin{definition}\label{def: sklt}
    A demi-normal pair $(X,\Delta_X)$ is called \emph{semi-Kawamata log terminal} (\emph{semi-klt}) if it is semi-dlt and $\lfloor\Delta_X\rfloor=0$.
\end{definition}

\begin{remark}\label{rmk: semiklt}
    Note that a semi-dlt pair $(X,\Delta_X)$ is semi-klt if and only if $(\Xbar,\Delta_{\Xbar})$ is klt; of course, the full normalization pair $(\Xbar,D_{\Xbar}+\Delta_{\Xbar})$ is dlt but not klt if $X$ is not normal. See~\S\ref{subsec: normal and demi-normal modification} for the notation.
\end{remark}

We conclude the subsection by briefly discussing the notion of $\QQ$-factoriality, as it is more subtle than in the algebraic case.

\begin{definition}[{\cite[Definition 2.38]{Fujino_MMPanalytic}, \cite[Definition 4.13]{Nakayama_ComplexPlurigenera}}]
    Let $X$ be a normal complex variety and let $K$ be a compact subset of $X$. Then $X$ is said to be $\QQ$-factorial at $K$ if every prime divisor defined on an open neighborhood $U$ of $K$ is $\QQ$-Cartier at any point $x\in K$. Let $\pi:X\rightarrow Y$ be a projective morphism and let $W$ be a compact subset of $Y$. If $X$ is $\QQ$-factorial at $\pi^{-1}(W)$, then we usually say that $X$ is $\QQ$-factorial over $W$.
\end{definition}

The main issue with $\QQ$-factoriality is that if $K'\subset K$ is an inclusion of compact subsets of $X$, and $X$ is $\QQ$-factorial at $K$, it is not automatically $\QQ$-factorial at $K'$, as there could be divisors defined on a neighborhood of $K'$ that cannot be extended to a neighborhood of $K$ (\cite[Remark 2.39]{Fujino_MMPanalytic}).

\subsubsection{Normalization and demi-normal modification}\label{subsec: normal and demi-normal modification}

Let $(X,\Delta_X)$ be a demi-normal pair, and $\pi_X\colon \Xbar\rightarrow X$ the normalization. The ideal sheaf of the conductor $D_X\subset X$ is the largest ideal sheaf of $\cO_X$ that is also an ideal sheaf of $\cO_{\Xbar}$, and it is equal to $\mathcal{H}om_X(\pi_*\cO_{\Xbar},\cO_X)$. The corresponding divisor in $\Xbar$ is also called conductor and denoted by $D_{\Xbar}$.
    
The \emph{normalization pair} is $(\Xbar,D_{\Xbar}+\Delta_{\Xbar})$ where $\Delta_{\Xbar}$ is the strict transform of $\Delta_X$.

We denote by $D_{\Xbar}^n$ the normalization of $D_{\Xbar}$. The map $D_{\Xbar}\rightarrow D_X$ between conductors is generically a double cover, and the associated bimeromorphic involution on $D_{\Xbar}$ extends to a regular involution $\mu\colon D_{\Xbar}^n\rightarrow D_{\Xbar}^n$. Moreover, $\mu$ is compatible with the different (\S\ref{subsec: adjunction}): $\mu^*\Diff_{D_{\Xbar}^n}(\Delta_{\Xbar})=\Diff_{D_{\Xbar}^n}(\Delta_{\Xbar})$, see~\cite[Proposition 5.12]{Kollar-singularities} for the algebraic case and \cite{Fujino_RelativeLogPluricanonical} for the analytic setting. The pair $(X,\Delta_X)$ is uniquely determined by the tuple $(\Xbar,D_{\Xbar},\Delta_{\Xbar},\mu)$ (see for instance \cite[Proposition 5.3]{Kollar-singularities}). In particular, two demi-normal pairs that are isomorphic in codimension 1 and have the same normalization are necessarily isomorphic.

Given a tuple $(\Xbar,D_{\Xbar}+\Delta_{\Xbar},\mu)$, deciding whether or not this arises as the normalization pair of a demi-normal pair is incredibly hard and not fully understood even for schemes. Koll\'ar's gluing theory \cite[\S5, \S9]{Kollar-singularities} was designed precisely to give a satisfactory answer to this problem in the context of algebraic slc pairs, and proved to be very successful. This complicated machinery is not yet available in the analytic world, hence we avoid using it.

The following construction is the demi-normal analogue of the normalization of a complex analytic space; see~\cite[Definition 5.1]{Kollar-singularities} for the algebraic case.

\begin{construction}[Demi-normal modification]\label{def: deminormal modification}
    Let $Z$ be a reduced complex analytic space. Suppose that there exists an open subspace $j\colon U\hookrightarrow Z$ with at most nodal singularities, and whose complement is analytic and of codimension at least 2. Then, $j_*\cO_U$ is a (locally) finite subsheaf of the coherent sheaf of $\cO_X$-modules $\pi_*\cO_{\Xbar}$, hence it is coherent (\cite[6.3]{Remmert_ComplexAnalyticSpaces}). To see the finite generation of $j_*\cO_U$, we can work locally and assume that both $\cO_X$ and $\cO_U$ are generated by their global sections, by Oka's Coherence Theorem~\cite[Theorem 0.14]{Fischer_BookAnalytic} and Cartan's A Theorem~\cite[Theorem 0.37]{Fischer_BookAnalytic}. Working around a Stein compact subset $K\subset X$ that satisfies property $\cP$ and using that $\cO_X(K)$ is Noetherian, we conclude that $j_*\cO_U(K)$ is also a finite $\cO_X(K)$-module, as wanted.
    Then, by construction, $\widetilde{Z}:=\Specan_Z j_*\cO_U$ is a demi-normal complex analytic space with a finite morphism $\delta_Z\colon \widetilde{Z}\rightarrow Z$, called the \emph{demi-normal modification} of $Z$.
\end{construction}

The following simple lemma shows that the demi-normal modification is independent of the open $U$, and satisfies a universal property.

\begin{lemma}\label{lem: universal property deminormal modification}
    With the notation of Definition~\ref{def: deminormal modification}, the following hold.
    \begin{enumerate}
        \item $\delta_Z$ is an isomorphism over $U$, and the normalizations of $\widetilde{Z}$ and $Z$ coincide.
        \item Let $f\colon X\rightarrow Z$ be a morphism from a demi-normal complex analytic space, and suppose that $f^{-1}(U)\subset X$ intersects all the irreducible components of the conductor of $X$. Then, $f$ factors uniquely through $\delta_Z$. Moreover, this property is satisfied if $f$ does not contract any component of the conductor. In particular, $\delta_Z\colon \widetilde{Z}\rightarrow Z$ is independent of $U$.
    \end{enumerate}
\end{lemma}
\begin{proof}
    The first point follows immediately from the definition. For the second point, note that giving a lift $X\rightarrow\widetilde{Z}$ of $f$ is equivalent to giving a factorization $\cO_Z\rightarrow j_*\cO_U\rightarrow f_*\cO_X$. Since $V=X\setminus D_X$ is normal, $f|_V\colon V\rightarrow Z$ factors through the normalization of $Z$, hence through $\widetilde{Z}$ by the first point of this lemma. Let $i\colon W=f^{-1}(U)\cup V\hookrightarrow X$, and note that its complement has codimension at least 2 in $X$. Since $X$ and $Z$ are Hausdorff, any factorization through $\widetilde{Z}$ is unique. In particular, $f|_W\colon W\rightarrow Z$ factors through $\delta_Z$. Since $X$ is demi-normal, we get the desired morphism $j_*\cO_U\rightarrow (f|_W)_*\cO_W=f_*i_*\cO_W\simeq f_*\cO_X$.
\end{proof}

\begin{example}\label{exm: deminormal modifications}
    We present two examples where the demi-normal modification is useful.
    \begin{enumerate}
    \item\label{exm: deminormal stein factorization} Let $X\rightarrow Y$ be a projective morphism from a demi-normal complex analytic space $X$ to a complex analytic space $Y$ that is nodal in codimension 1, and let $X\rightarrow Z\rightarrow Y$ be the Stein factorization. Suppose that $X\rightarrow Y$ has connected fibers over a Zariski-open subset $U\subset Y$ whose complement has codimension at least 2, and $f$ does not contract any component of the conductor. Then, $Z\rightarrow Y$ is an isomorphism over $U$, hence $Z$ is also nodal in codimension 1. Since $X$ is demi-normal, $X\rightarrow Z$ factors through the demi-normal modification $\widetilde{Z}\rightarrow Z$, which is a finite map. From the properties of the Stein factorization, we get that $\widetilde{Z}\simeq Z$, that is, $Z$ is demi-normal.
    \item\label{exm: deminormal boundary} Let $X$ be a normal complex analytic space and $S\subset X$ a reduced subspace of pure codimension 1 such that $(X,S+B)$ is an lc pair for some $\QQ$-divisor $B$. Then, $S$ is nodal in codimension 1 by~\cite[Proposition 16.6]{K+92}, hence we can take the universal demi-normal modification $\widetilde{S}\rightarrow S$. In general $S$ does not satisfy Serre's condition $S_2$, even when $S$ is realized as the conductor of a normalization pair. See~\cite[Example 17.5.2]{K+92} or the next example for instances where the $S_2$ property fails.
    \item\label{exm: conductor is not always S2}
    Consider the product $S=E\times\PP^1$ where $E$ is a smooth curve of genus 1. Let $L$ be an ample line bundle on $S$, and let $X$ be the cone over $S$ corresponding to $L$. Its natural resolution is the total space of $L$ itself. Let $D=E\times\{0,\infty\}$, and $\Delta_X$ be the cone over $D$ corresponding to $L|_D$. Since $K_S+D\sim 0$, the pair $(X,\Delta_X)$ is an lc pair that is not dlt by~\cite[Lemma 3.1]{Kollar-singularities}. In fact, $\Delta_X=\lfloor\Delta_X\rfloor$ is the union of two surfaces meeting at one point, hence not $S_2$.
    \end{enumerate}
\end{example}

\begin{remark}\label{rmk: double cover}
    There exists another useful construction involving demi-normal spaces. Given a demi-normal complex analytic space $X$, there exists a double-cover $g\colon Y\rightarrow X$ that is étale in codimension 1 and such that $Y$ is simple normal crossing in codimension 1. In particular, the irreducible components of $Y$ are normal in codimension 1. See~\cite[5.23]{Kollar-singularities} for the construction, which works in the complex setting without changes.
\end{remark}

In the proof of Theorems~\ref{thm: general version} and~\ref{thm: equivariant sklt}, we will need to extend a meromorphic involution to a holomorphic involution. The following analytic analogue of~\cite[Theorem 11.40]{Kollar-moduli} will be crucial. Even though the proof is almost identical to the algebraic case, we report it here for the sake of completeness.

\begin{theorem}[{\cite[Theorem 11.40]{Kollar-moduli}}]\label{thm: extending involution}
    Let $T$ be a complex analytic space, and let $f_i\colon(X_i,\Delta_{X_i})\rightarrow T$ be projective morphisms from lc pairs such that $K_{X_i}+\Delta_{X_i}$ is $f_i$-ample. Let $Z_T\subset T$ be a closed subset and set $Z_i:=f_i^{-1}(Z_T)$. Let $\mu^0:(X_1\setminus Z_1,\Delta_{X_1}|_{X_1\setminus Z_1})\rightarrow(X_2\setminus Z_2,\Delta_{X_2}|_{X_2\setminus Z_2})$ be an isomorphism. Assume that none of the log centers of $(X_i,\Delta_{X_i})$ is contained in $Z_i$ for $i=1,2$. Then $\mu^0$ extends to an isomorphism $\mu:(X_1,\Delta_{X_1})\rightarrow (X_2,\Delta_{X_2})$.
\end{theorem}
\begin{proof}
    By the assumption on log centers, $Z_i$ does not contain any component of $X_i$ or $\Delta_i$. Let $X_1\xleftarrow{p_1}\Gamma\xrightarrow{p_2}X_2$ be the normalization of the closure of the graph of $\mu^0$. Let $\Delta_{\Gamma}:=p_{1*}^{-1}\Delta_1=p_{2*}^{-1}\Delta_2$, and write
    \[
        K_{\Gamma}+\Delta_{\Gamma}\sim_{\QQ}p_i^*(K_{X_i}+\Delta_{X_i})+E_i,
    \]
    where $E_i$ is $p_i$-exceptional and $p_i(\Supp E_i)\subset Z_i$ by construction. Note that $E_i$ is effective by the assumption on the log centers. By the ampleness assumption and the Negativity Lemma~\cite[Corollary 2.16]{Enokizono-Hashizume-analytic-MMP} applied to $p_1$ and $p_2$ separately, we get $E_1=E_2$. Since $\Gamma\rightarrow X_1\times_T X_2$ is finite, any curve contracted by $p_2:\Gamma\rightarrow X_2$ is not contracted by $p_1:\Gamma\rightarrow X_1$. Since $p_1^*(K_{X_1}+\Delta_{X_1})=p_2^*(K_{X_2}+\Delta_{X_2})$ and both log canonical divisors are ample, it follows that $p_2$ cannot contract any curve, and the same holds for $p_1$. Therefore, both $p_1$ and $p_2$ are finite bimeromorphic morphisms between normal complex analytic spaces, hence isomorphisms by Zariski's Main Theorem. It follows that $\mu=p_2\circ p_1^{-1}$ is an isomorphism extending $\mu^0$.
\end{proof}

\subsubsection{Adjunction}\label{subsec: adjunction}

Adjunction for complex analytic spaces works in the same way as for schemes; see for instance~\cite{Kollar_Mori_1998,Fujino_MMPanalytic,Das_Hacon_Pau_4dimensionKahler,Fujino_inversion_adjunction_analytic}. In this subsection we limit ourselves to collecting the following results on adjunction for (semi-)dlt pairs and recalling the description of the log canonical ring of an slc pair in terms of the log canonical ring of its normalization and the gluing involution on the conductor (\S\ref{subsec: normal and demi-normal modification}).

\begin{lemma}\label{lem: semidlt adjunction}
    Let $(X,\Delta_X)$ be a dlt pair, $S\leq\lfloor\Delta_X\rfloor$, and $B=\Delta_X-S$. Then,
    \begin{enumerate}
        \item $S$ is demi-normal, and $(S,\Diff_{S}(B))$ is semi-dlt.
        \item If $S=\lfloor\Delta_X\rfloor$ (equivalently, $\lfloor B\rfloor=0$), then $(S,\Diff_{S}(B))$ is semi-klt.
    \end{enumerate}
    In particular, if $(Z,\Delta_Z)$ is a semi-dlt pair and $(\Zbar,D_{\Zbar}+\Delta_{\Zbar})$ is the normalization pair, then $(D_{\Zbar},\Diff_{D_{\Zbar}}(\Delta_{\Zbar}))$ is semi-dlt. Moreover, $(D_{\Zbar},\Diff_{D_{\Zbar}}(\Delta_{\Zbar}))$ is semi-klt if and only if $(Z,\Delta_Z)$ is semi-klt along its conductor $D_Z$.
\end{lemma}
\begin{proof}
    The first point is well-known, see for instance~\cite[Definition 5.19]{Kollar-singularities} for the algebraic case.
    
    For the second point, we can reduce to the case where $X$ is a surface. In the algebraic case, this is achieved by localizing at codimension 1 points of $S$. In the complex analytic case, we can embed analytically-locally $X$ into a complex manifold, and intersect with a very general transverse codimension 2 subspace. We can then assume $(X,S+B)$ to be a surface dlt pair, hence analytically-locally $\QQ$-factorial (see for instance \cite[Proposition 4.11]{Kollar_Mori_1998}). By $\QQ$-factoriality, we have that $\Diff_{S}(B)=\Diff_{S}(0)+B|_S$, and by the first point it is enough to show that $\lfloor\Diff_S(B)\rfloor=0$. By~\cite[Theorem 3.36]{Kollar-singularities}, if $p\in S$ is a point at which $\Diff_S(0)$ has coefficient 1, then $S$ is smooth at $p$. Then, $(X,S+B)$ is plt locally around $p$ by~\cite[Proposition 2.42]{Kollar_Mori_1998}, hence $(S^n,\Diff_{S^n}(B))$ is klt around $p$ by adjunction (\cite[Theorem 5.50]{Kollar_Mori_1998}); in particular, the coefficient of $\Diff_{S}(0)$ is actually strictly smaller than 1. Now, let $p$ be in the support of $B|_S$. Since we are working with surfaces, the characterization of lc centers on dlt pairs implies that $p$ is again a smooth point of $S$. It follows again by adjunction (\cite[Theorem 5.50]{Kollar_Mori_1998}) that the boundary pair is klt.

    Finally, let $(Z,\Delta_Z)$ be a semi-dlt pair. Then, its normalization pair is dlt, and $(\Zbar,\Delta_{\Zbar})$ is klt if and only if $(Z,\Delta_Z)$ is semi-klt (Remark~\ref{rmk: semiklt}). Then, the result follows from the first two points of the statement.
\end{proof}

See~\cite[Proposition 5.8]{Kollar-singularities} for a proof of the following lemma in the case where $T$ is a reduced point; see also~\cite{Fujino_RelativeLogPluricanonical}.

\begin{lemma}\label{lem: log canonical ring slc pair}
    Let $(X,\Delta_X)$ be an slc pair with normalization pair $(\Xbar,D_{\Xbar}+\Delta_{\Xbar})$, where $D_{\Xbar}$ is the conductor, and let $f\colon X\rightarrow T$ be a projective morphism to a complex analytic space $T$. Let $\mu$ be the gluing involution on $(D_{\Xbar}^n,\Diff_{D_{\Xbar}^n}(\Delta_{\Xbar}))$, and let $R_{D_{\Xbar}^n}:=R(D_{\Xbar}^n/T,K_{D_{\Xbar}^n}+\Diff_{D_{\Xbar}^n}(\Delta_{\Xbar}))$. Let $\mu$ act on $R_{D_{\Xbar}^n}$ via $\mu\cdot\alpha=(-1)^m\mu^*\alpha$ for all $\alpha\in R_{D_{\Xbar}^n}$ of degree $m$. Then, the inclusion
    \[
    \begin{tikzcd}
        R_X:=R(X/T,K_X+\Delta_X)\arrow[r,hookrightarrow] & R(\Xbar/T,K_{\Xbar}+D_{\Xbar}+\Delta_{\Xbar})=:R_{\Xbar}
    \end{tikzcd}
    \]
    yields an identification
    \[
        R_X=\res_n^{-1}(R_{D_{\Xbar}^n}^\mu),
    \]
    where $R_{D_{\Xbar}^n}^{\mu}$ denotes the subalgebra of $\mu$-invariants, and
    \[
    \begin{tikzcd}
        \res_n\colon R_{\Xbar}\arrow[r] & R_{D_{\Xbar}^n}
    \end{tikzcd}
    \]
    is the natural residue (or restriction) map.
\end{lemma}

\begin{remark}
    The residue map $\res_n$ factors through $\res\colon R_{\Xbar}\rightarrow R_{\widetilde{D}_{\Xbar}}$, where $\widetilde{D}_{\Xbar}$ is the demi-normal modification of $D_{\Xbar}$ (see Definition~\ref{def: deminormal modification} and Example~\ref{exm: deminormal modifications}\eqref{exm: deminormal boundary}) and $R_{\widetilde{D}_{\Xbar}}$ the corresponding log canonical $\cO_T$-algebra. In particular, if $(X,\Delta_X)$ is semi-dlt then $\res_n$ factors through $\res\colon R_{\Xbar}\rightarrow R_{D_{\Xbar}}$.
\end{remark}

\subsection{MMP for projective morphisms between complex analytic spaces}\label{subsec: basics MMP}

In this subsection we briefly recall the main notions and results of the minimal model program for projective morphisms between complex analytic spaces.
We limit ourselves to discuss the principal differences with the algebraic case and provide references to the various statements. The main technical point is that most definitions and constructions are performed around compact subsets in order to have a well-behaved notion. We refer to~\cite{Fujino_MMPanalytic} and~\cite[\S3]{Enokizono-Hashizume-analytic-MMP} for details.
We also extend some results on the existence of minimal models to the analytic setting, such as~\cite[Theorem 1.1]{Hacon_Xu_dslt} in Theorem~\ref{thm: Hacon-Xu analytic}; see also Lemmas~\ref{lem: termination very exceptional MMP} and~\ref{lem: existence minimal model if log canonical}.

First of all, it is still possible to define relative Kleiman-Mori cones, but in order to work with finite dimensional spaces we need to work around particularly nice compact subsets.

\begin{construction}[{\cite[\S4]{Fujino_MMPanalytic}, \cite[page 70]{Nakayama_ZariskiDecomposition}}]
    Let $f\colon X\rightarrow T$ be a projective morphism of complex analytic spaces, and $W\subset T$ a compact subset satisfying property $\cP$. Let $Z_1(X/T;W)$ be the free abelian group generated by curves in $X$ that are contracted to a point of $W$ via $f$. Every line bundle defined over $f^{-1}(U)$ for some open $U\supset W$ can be paired with elements in $Z_1(X/T;W)$, and defines a numerical equivalence relation on $\mathrm{Pic}(f^{-1}(U))$. Taking the quotient by this equivalence relation, and the colimit over all open neighborhoods $U$ of $W$, we get an abelian group $A^1(X/T;W)$, which is finitely generated by~\cite[II 5.19 Lemma]{Nakayama_ZariskiDecomposition} and~\cite[Lemma 4.1]{Fujino_MMPanalytic}. We define $A_1(X/T;W)$ as the image of $Z_1(X/T;W)\rightarrow\Hom(A^1(X/T;W),\ZZ)$, and $N_1(X/T;W):=A_1(X/T;W)\otimes\RR$, $N^1(X/T;W):=A^1(X/T;W)\otimes\RR$. Finally, we define the \emph{Kleiman-Mori cone} $\overline{\mathrm{NE}}(X/T;W)$ of $f$ over $W$ as the closure of the convex cone in $N_1(X/T;W)$ spanned by the classes of integral curves contracted to a point of $W$. We remark that in this construction the assumption of $W$ being Stein compact is not needed.
\end{construction}

\subsubsection{MMP steps}

In the complex analytic setting, \emph{MMP steps} (with scaling) for lc pairs are defined as in the algebraic case, with the difference that one defines them only on a neighborhood of a compact subsets satisfying $\cP$, see~\cite[Definition 3.5]{Enokizono-Hashizume-analytic-MMP}. Since we will be interested in the minimal model program for slc pairs, we partially extend the definition to include particularly nice MMP steps between slc pairs.

\begin{definition}[{\cite[Definition 2.6]{Hashizume_sdlt_models}, \cite[Definition 3.5]{Enokizono-Hashizume-analytic-MMP}}]\label{def:slc-MMP steps}
    Let $(X,\Delta_X)$ be an slc pair, and let $f\colon X\rightarrow T$ be a projective morphism with $T$ a complex analytic space and $W\subset T$ a compact subset satisfying $\cP$. Then, an \emph{MMP step over $T$ around $W$} for $(X,\Delta_X)$ is a commutative diagram
    \begin{equation}\label{diag: slc MMP step}
    \begin{tikzcd}
        (X,\Delta_X)\arrow[rr,dashed,"\psi"]\arrow[rd,"\varphi"'] & & (Y,\Delta_Y)\arrow[ld,"\varphi'"]\\
        & Z
    \end{tikzcd}
    \end{equation}
    projective over $T$, such that
    \begin{enumerate}
        \item\label{slc_step:1} $(Y,\Delta_Y)$ is an slc pair;
        \item\label{slc_step:2} $\psi$ is bimeromorphic and $\Ex(\psi)$ does not contain any component of the conductor;
        \item $\Delta_Y=\psi_*\Delta_X$;
        \item the morphisms $Z\rightarrow T$ and $Y\rightarrow T$ are projective, $\varphi$ and $\varphi'$ are generically finite, and $\varphi'$ has no exceptional divisors;
        \item $-(K_X+\Delta_X)$ and $(K_Y+\Delta_Y)$ are ample over $Z$.
    \end{enumerate}
    Let $H$ be a $\QQ$-Cartier divisor on $X$. An \emph{MMP step over $T$ around $W$ with scaling of $H$} is an MMP step as in~\eqref{diag: slc MMP step} satisfying the following two conditions:
    \begin{enumerate}[start=6]
        \item $K_X+\Delta_X+tH$ is nef for some $t\in\QQ_{>0}$, and
        \item if we put
        \[
            \lambda:=\mathrm{inf}\{t\in\RR_{\geq0}: K_X+\Delta_X+tH\text{ is nef over }W\}
        \]
        then $(K_X+\Delta_X+\lambda H)\cdot C=0$ for any curve $C\subset\pi^{-1}W$ contracted by $f$. 
    \end{enumerate}
    A \emph{sequence of MMP steps over $T$ around $W$ with scalings of $H$} is a sequence
    \[
        (X,\Delta_X)\dashrightarrow(X_1,\Delta_{X_1})\dashrightarrow\ldots\dashrightarrow(X_i,\Delta_{X_i})\dashrightarrow\ldots
    \]
    possibly defined only after restricting $T$ around $W$ at each step, such that $(X_i,\Delta_{X_i})\dashrightarrow(X_{i+1},\Delta_{X_{i+1}})$ is an MMP step over $T$ around $W$ with scalings of $H_i$, where $H_i=\phi_{i*}H$ and $\phi_i:X\rightarrow X_i$ is the composite of the first $i$ steps.
\end{definition}

\begin{remark}\label{rmk: normalization MMP step}
    It can be shown as in~\cite[Lemma 12]{Ambro-Kollar-slc-MMP} that a diagram as~\eqref{diag: slc MMP step} satisfying properties~\eqref{slc_step:1} and~\eqref{slc_step:2} is an MMP step over $T$ for $(X,\Delta_X)$ if and only if the diagram obtained by normalization is an MMP step over $T$ for $(\Xbar,D_{\Xbar}+\Delta_{\Xbar})$. In the literature, $\varphi$ is often required to be the contraction of a $(K_X+\Delta_X)$-negative extremal ray of $\overline{\mathrm{NE}}(X/Y;W)$, as in~\cite[Definition 3.5]{Enokizono-Hashizume-analytic-MMP}. Moreover, one usually further assumes $T$ to be Stein and $W$ to be Stein compact (satisfying $\cP$), as this is the setup in which the existence of MMP steps is guaranteed for log canonical pairs. We will need a slightly more general setting, particularly in Theorem~\ref{thm: semidlt modifications} and Lemma~\ref{lem: termination very exceptional MMP}.
\end{remark}

\begin{remark}\label{rmk: rationality nef threshold}
    We will only use Definition~\ref{def:slc-MMP steps} when $H$ is $f$-ample. Then, the rationality theorem (\cite[Theorem 7.1]{Fujino_MMPanalytic}, \cite[Theorem 9.1]{Fujino_QuasiLogAnalytic}) implies that the \emph{nef threshold} $\lambda$ is a rational number. In particular, we will not have to worry about defining $\RR$-divisors in this paper.
\end{remark}

\subsubsection{Minimal models}

In the complex analytic setting, \emph{weak log canonical models} and \emph{minimal models} are defined in the same way as for schemes, with the only difference being that to have a well-behaved notion one restricts to working around compact subsets satisfying property $\cP$, often assumed to be Stein compact; see~\cite[Definition 11.4]{Fujino_MMPanalytic} or~\cite[\S3.1]{Enokizono-Hashizume-analytic-MMP}. Since the \emph{log canonical model} (or ample model, or stable model) is unique when it exists (\cite[Lemma 11.3]{Fujino_MMPanalytic}), the construction of a relative log canonical model of an lc pair $(X,\Delta_X)$ with a projective morphism $f\colon X\rightarrow T$ is local on $T$. In particular, we can always work around compact subsets of $T$ and then patch the various log canonical models together using uniqueness, hence there is no difference between the algebraic and analytic setting. As in the algebraic case, if $X\dashrightarrow X^c$ is the log canonical model of an lc pair $(X,\Delta_X)$ projective over $T$, then $R_X:=R(X/T,K_X+\Delta_X)$ is a locally finitely presented $\cO_T$-algebra and $X^c\simeq\Projan_T R(X/T,K_X+\Delta_X)$. Conversely, if $R_X$ is locally finitely generated (hence locally finitely presented by Lemma~\ref{lem: finite presentation}), its analytic Proj is the log canonical model. See Proposition~\ref{prop: Y=proj RX good stable model} and Theorem~\ref{thm: projRX is stable model} for a generalization to slc pairs.

The main tools of the minimal model program for log canonical pairs have been translated to the complex analytic realm. For instance, the main results of~\cite{BCHM} (\cite{Fujino_MMPanalytic}), vanishing theorems (\cite{Fujino_MMPanalytic}), the cone and contraction theorem (\cite{Fujino_ConeContractionAnalytic}), existence of ($\QQ$-factorial) dlt modifications (\cite[Theorems 1.21, 1.24, 1.27]{Fujino_MMPanalytic}), and existence of flips (\cite[Theorem 1.2]{Enokizono-Hashizume-analytic-MMP}, \cite[Theorem 1.7]{Fujino_RelativeLogPluricanonical}) are all available in the analytic setting as well.

For slc pairs, vanishing theorems, the cone and contraction theorem are also available, as a consequence of the theory of quasi-log structures on complex analytic spaces~\cite{Fujino_QuasiLogAnalytic}. Moreover, as a consequence of the finiteness of B-representations (\cite[Theorem 1.2]{Fujino_RelativeLogPluricanonical}), for any slc pair $(X,\Delta_X)$ with a projective morphism $f\colon X\rightarrow T$ the log canonical divisor class $K_X+\Delta_X$ is $f$-semi-ample if and only if the log canonical divisor class $K_{\Xbar}+D_{\Xbar}+\Delta_{\Xbar}$ of the normalization pair $(\Xbar,D_{\Xbar}+\Delta_{\Xbar})$ is $f$-semi-ample (\cite[Theorem 1.1]{Fujino_RelativeLogPluricanonical}).

We will need the following lemmas, that guarantee the existence of good minimal models in some cases. Recall that a given a projective morphism $\varphi:X\rightarrow Y$ between normal analytic varieties with $\varphi_*\cO_X\simeq\cO_Y$, a $\QQ$-divisor $D$ on $X$ is \emph{very
exceptional} over $Y$ if $D$ is vertical and for any prime divisor $Q$ on $Y$, there is a prime
divisor $P$ on $X$ such that $\varphi(P)=Q$ and $P$ is not a component of $D$. This notion is well-behaved even in the complex setting, and coincides with the notion of exceptional divisor when $\varphi$ is bimeromorphic; see~\cite[\S2.4]{Enokizono-Hashizume-analytic-MMP}.

\begin{lemma}[{\cite[Theorem 3.4]{Birkar_existenceflips}}]\label{lem: termination very exceptional MMP}
    Let $(X,\Delta_X)$ be a dlt pair, $f:X\rightarrow T$ a projective morphism to a Stein space, and $W\subset T$ a Stein compact subset satisfying $\cP$ such that $X$ is $\QQ$-factorial over $W$. Let $\varphi:X\rightarrow Y$ and $g:Y\rightarrow T$ be projective morphisms with $Y$ a complex analytic space and $f=g\circ\varphi$. Suppose that $K_X+\Delta_X\sim_{\QQ,Y}E$ for some effective and very exceptional divisor $E$ over $Y$. Then, any sequence of MMP steps over $Y$ around $g^{-1}(W)$ with scaling of a general $\varphi$-ample $\QQ$-divisor terminates with a minimal model $(Z,\Delta_Z)$ over $Y$ around $g^{-1}(W)$ that contracts all components of $E$, and such that $K_Z+\Delta_Z\sim_{\QQ,Y}0$.
\end{lemma}
\begin{proof}
    Let $H$ be a general $\varphi$-ample $\QQ$-divisor on $X$. The first step would be to run an MPP over $Y$ around $g^{-1}(W)$ with scalings of $H$. A technical difficulty is that $g^{-1}(W)$ is not Stein compact, hence we cannot directly apply the results of the minimal model program for projective morphisms between complex analytic spaces. To solve this problem, we use an idea of Fujino~\cite[proof of Theorem 1.27]{Fujino_MMPanalytic}. Let $A$ be a general $g$-ample divisor on $Y$ such that $A\cdot C>2\ \!\dim\ \! Y$ for every curve $C$ that is contracted by $g$; this exists after suitably restricting $T$ around $W$. Note that $\varphi^*A+H$ is $f$-ample. By~\cite[Lemma 9.4]{Fujino_MMPanalytic}, a $(K_X+\Delta_X+\varphi^*A)$-MMP step over $T$ with scalings of $H$ is also a $(K_X+\Delta_X)$-MMP step over $Y$ with scalings of $\varphi^*A+H\sim_{\QQ,Y}H$. After restricting $T$ around $W$ at each step, we can construct a sequence of $(K_X+\Delta_X+\varphi^*A)$-MMP steps over $T$ around $W$ with scalings of $H$, where all pairs are dlt and $\QQ$-factorial over $W$ at each step. By the discussion above, we can see this sequence as an MMP over $Y$ around $g^{-1}(W)$ with scalings of $H$. Recall that $K_X+\Delta_X$ is $\QQ$-linearly equivalent over $Y$ to the very exceptional divisor $E$, and this holds at each step of the MMP. Then, the same argument as~\cite[Theorem 3.4]{Birkar_existenceflips} implies that the MMP terminates with a dlt minimal model $(Z,\Delta_Z)$ over $Y$, and contracts all components of $E$. We recall the reasoning of~\cite[Theorem 3.4]{Birkar_existenceflips} here for the convenience of the reader. Suppose by contradiction that the MMP sequence at the normalization level is infinite. We can assume it to consists only of flips, and by the analytic analogs~\cite{Fujino_MMPanalytic} of~\cite{BCHM}, the limit of the nef thresholds is 0. Let $\phi_i:X\rightarrow X_i$ be the composition of the first $i$ steps of the MMP, and let $H_i:=\phi_{i*}H$. Since we are assuming all steps to be flips, the strict transform of $K_{X_i}+\Delta_{X_i}+\lambda_i H_i$ along $\phi_i^{-1}$ is equal to $K_X+\Delta_X+\lambda_i H$. Since $K_{X_i}+\Delta_{X_i}+\lambda_i H_i$ is nef on $X_i$ and the $\lambda_i$'s converge to 0, we have that $K_X+\Delta_X$ is a limit of movable divisors. It follows that for any prime divisor $S$ and very general curve $C\subset S$, we have $(K_X+\Delta_X)\cdot C\geq0$, which however would contradict the negativity lemma if $E\not=0$ (\cite[Lemma 3.3]{Birkar_existenceflips}, \cite[Lemma 2.14]{Enokizono-Hashizume-analytic-MMP}). This shows that the MMP terminates with a good minimal model and contracts all components of $E$. In particular, $K_Z+\Delta_Z\sim_{\QQ,Y}0$.
\end{proof}

\begin{lemma}\label{lem: existence minimal model if log canonical}
    Let $f\colon X\rightarrow T$ be a projective morphism of complex analytic spaces and $(X,\Delta_X)$ an lc pair admitting a log canonical model $\phi\colon (X,\Delta_X)\dashrightarrow (Y,\Delta_Y)$ over $T$. Suppose $T$ to be Stein and $W\subset T$ a Stein compact subset satisfying property $\cP$. Then, $(X,\Delta_X)$ admits a $\QQ$-factorial dlt good minimal model over $T$ around $W$. Moreover, if $(X,\Delta_X)$ is dlt and $\QQ$-factorial over $W$, then we can reach a good minimal model by performing a finite sequence of MMP steps over $T$ around $W$ with scalings of an $f$-ample divisor.
\end{lemma}
\begin{proof}
    Since the statement is local around $W$, we will tacitly restrict $T$ around it when needed in order to take log resolutions and apply results of the minimal model program; in particular, we will assume $T$ to be Stein. Let $p:X'\rightarrow X$ be a log resolution that resolves the meromorphic map $\phi$ to the log canonical model, and let $\Delta_{X'}$ be the boundary satisfying $K_{X'}+\Delta_{X'}=p^*(K_X+\Delta_X)+E$ for some $p$-excepional divisor $E\geq0$ having no components in common with $\Delta_{X'}$. Then, the relative log canonical algebras of $(X',\Delta_{X'})$ and $(X,\Delta_X)$ coincide, $(Y,\Delta_Y)$ is the log canonical model of $(X',\Delta_{X'})$, and a good minimal model of $(X',\Delta_{X'})$ would be a good minimal model for $(X,\Delta_X)$ as well. Thus, from now on, we can assume $\phi$ to be a morphism and $(X,\Delta_X)$ to be dlt and $\QQ$-factorial over $W$. By the Negativity Lemma~\cite[Corollary 2.16]{Enokizono-Hashizume-analytic-MMP}, we have that $K_X+\Delta_X=\phi^*(K_Y+\Delta_Y)+E$ for some effective $\phi$-exceptional divisor $E$. By Lemma~\ref{lem: termination very exceptional MMP} and its proof, we can run an MMP over $Y$ around $f^{-1}(W)$ with scalings of an ample divisor, that is also an MMP over $T$ around $W$, which terminates with a minimal model over $Y$. Note that this is also a good minimal model over $T$. In the algebraic case, the existence of a good minimal model implies the termination of every MMP with scaling of a relatively ample divisor, see~\cite[Corollary 2.9]{Hacon_Xu_dslt}. The same proof with the techniques of~\cite{Fujino_MMPanalytic} works in the analytic setting.
\end{proof}

We conclude by stating a generalization of~\cite[Theorem 1.1]{Hacon_Xu_dslt} to the complex analytic setting. We only sketch the proof, as it follows from known arguments.

\begin{theorem}[{\cite[Theorem 1.1]{Hacon_Xu_dslt}, \cite[Theorem 5.8]{Enokizono-Hashizume-analytic-MMP}}]\label{thm: Hacon-Xu analytic}
    Let $f\colon X\rightarrow T$ be a projective morphism from a normal complex analytic space $X$ to a Stein space $T$, and let $W\subset T$ be a Stein compact subset satisfying $\cP$. Let $(X,\Delta_X)$ be a dlt pair and $S=\lfloor\Delta_X\rfloor$ the non-klt locus. Assume that there exists a Zariski-open subset $T^0\subset T$ such that $(X^0,\Delta_{X^0})=(X,\Delta_X)\times_TT^0$ has a good minimal model over $T^0$, and that any stratum of $S$ intersects $X^0$. Let $H$ be an $f$-ample $\QQ$-divisor on $X$ such that $(X,\Delta_X+H)$ is lc and $K_X+\Delta_X+H$ is nef over $W$. Then, after possibly restricting $T$ around $W$ further, there exists a finite sequence of steps of a $(K_X+\Delta_X)$-MMP over $Y$ around $W$ with scaling of $H$
    \[
        (X,\Delta_X)\dashrightarrow(X_1,\Delta_{X_1})\dashrightarrow\ldots\dashrightarrow(X_m,\Delta_{X_m})
    \]
    where the lc pair $(X_m,B_m)$ is a good log minimal model of $(X,\Delta_X)$ over $T$ around $W$.
\end{theorem}
\begin{proof}
    This follows from~\cite[Theorem 5.8]{Enokizono-Hashizume-analytic-MMP} and the argument in~\cite[Remark 3.7]{Hashizume_ScalingFinitenessLogAbundant}. Indeed, after running a finite number of steps of the MMP with scalings and restricting $T$ around $W$, we can assume that $K_{X^0}+\Delta_{X^0}$ is semi-ample over $T^0$. Then, by the assumption on the strata of $S$, $K_X+\Delta_X$ is log abundant (\cite[Definition 2.16]{Hashizume_ScalingFinitenessLogAbundant}) and it stays so after any further step of the MMP that is an isomorphism over $T^0$. Therefore, the MMP with scalings terminates by~\cite[Theorem 5.8]{Enokizono-Hashizume-analytic-MMP}.
\end{proof}

\subsubsection{Kawamata-Viehweg vanishing}

An important ingredient in the proof of Theorem~\ref{thm: equivariant sklt} and for the MMP in general is Kawamata-Viehweg vanishing. Fortunately, this is available in the complex analytic setting as well; we record the statement below, see also~\cite{Fujino_VanishingAnalytic}.

\begin{theorem}[{\cite[Theorem 2.21]{Das_Hacon_Kahler3foldsMMP}}]
    Let $f\colon X\rightarrow Y$ be a projective surjective morphism of normal complex analytic spaces, with $(X,\Delta_X)$ a klt pair. Let $D$ be an $\QQ$-Cartier integral divisor on $X$, such that $D-(K_X+\Delta_X)$ is $f$-nef and $f$-big. Then,
    \[
        R^if_*\cO_X(D)=0\qquad\text{for all $i>0$.}
    \]
\end{theorem}

As an immediate but crucial consequence, we get the following lemma.

\begin{lemma}\label{lem: surjectivity restriction minimal model}
    Let $f\colon Y\rightarrow T$ be a projective morphism between complex analytic spaces, and let $(Y,S+B)$ be a $\QQ$-factorial dlt pair with $\lfloor S+B\rfloor=S$. Suppose that $K_Y+S+B$ is $f$-nef and $f$-big. Then, for $\ell$ sufficiently divisible, $\ell(K_Y+S+B)$ is Cartier and the restriction map
    \[
    \begin{tikzcd}
        f_*\cO_Y(\ell(K_Y+S+B))\arrow[r,twoheadrightarrow] & f_*\cO_{S}(\ell(K_S+\Diff_S(B)))
    \end{tikzcd}
    \]
    is surjective.
\end{lemma}

\subsection{B-bimeromorphic maps}\label{subsec: B-bimeromorphic}

We recall the basic notion of \emph{B-bimeromorphic maps}; here, the `B' stands for `boundary', as they can be thought as bimeromorphic maps preserving the boundary. This is the complex analytic analogue of the notion of B-birational maps in the algebraic setting, which has proved to be of crucial importance when dealing with slc pairs~\cite{FujinoAbundanceSLC,FG_fin_Brepr,Hacon_Xu_finBrepr_slc_abundance,Fujino_RelativeLogPluricanonical}.

\begin{definition}[Relative B-bimeromorphic maps, {\cite[Definition 2.17]{Fujino_RelativeLogPluricanonical}}]\label{def: B-birational maps}
    Let $f\colon X\rightarrow T$ and $g\colon Y\rightarrow T$ be projective morphisms to a complex analytic space $T$, and let $(X,\Delta_X)$ and $(Y,\Delta_Y)$ be log canonical pairs. We say that a bimeromorphic map $\sigma\colon (X,\Delta_X)\dashrightarrow(Y,\Delta_Y)$ over $T$ is \emph{B-bimeromorphic} if there exists a commutative diagram
    \[
    \begin{tikzcd}
        & Z\arrow[ld,"\alpha"']\arrow[rd,"\beta"]\\
        X\arrow[rr,dashed,"\sigma"] & & Y
    \end{tikzcd}
    \]
    where $Z$ is normal, $\alpha$, $\beta$ are proper bimeromorphic morphisms, and
    \[
        \alpha^*(K_X+\Delta_X)=\beta^*(K_Y+\Delta_Y).
    \]
    We denote by $\Bim(X/T,\Delta_X)$ the group of B-bimeromorphic self-maps of $(X,\Delta_X)$ over $T$, where the group operation is the composition.
\end{definition}


\begin{remark}\label{rmk: Bir is a crepant invariant}
    An important property of $\Bim(X,\Delta_X)$ is that it is a crepant bimeromorphic invariant, in the following sense. Let $\psi\colon (X,\Delta_X)\rightarrow(Y,\Delta_Y)$ be a proper crepant bimeromorphic morphism between lc pairs that are projective over a complex analyic space $T$, then $\Bim(Y/T,\Delta_Y)\simeq\Bim(X/T,\Delta_X)$ via $g\mapsto\psi^{-1}\circ g\circ\psi$; see~\cite[Remark 2.12]{FG_fin_Brepr}.
    We will apply this observation when $(X,\Delta_X)$ is a minimal model over $T$, and $(Y,\Delta_Y)$ is its relative log canonical model.
\end{remark}

Another property of B-bimeromorphic self-maps we will use is that they act on the spaces of pluri-canonical sections, thus on a suitable Veronese subalgebra of the log canonical ring. This follows immediately from the definition of B-bimeromorphic self-maps over $T$, hence we omit the proof of the following lemma; see for instance~\cite[Definitions 2.17, 2.22]{Fujino_RelativeLogPluricanonical} or~\cite[Definition 2.14]{FG_fin_Brepr}.

\begin{lemma}\label{lem: Bir acts on log canonical}
    Let $X$ be a normal complex analytic space such that $(X,\Delta_X)$ is log canonical, and $f\colon X\rightarrow T$ a projective morphism to a complex analytic space. Let $m$ be a non-negative integer such that $m(K_X+\Delta_X)$ is Cartier. Then, any B-bimeromorphic map $\sigma\in\Bim(X/T,\Delta_X)$ over $T$ naturally induces a linear automorphism of $f_*\cO_{X}(m(K_X+\Delta_X))$.
    
    In particular, if $\ell$ is the smallest positive integer such that $\ell(K_X+\Delta_X)$ is Cartier, then $\Bim(X/T,\Delta_X)$ acts on $R(X/T,K_X+\Delta_X)^{(\ell)}$.
\end{lemma}

Although we will not use this, it is worth mentioning that the induced homomorphisms $\rho_m:\Bim(X/T,\Delta_X)\rightarrow\Aut_{\cO_T}(f_*\cO_X(m(K_X+\Delta_X)))$ are called \emph{B-pluricanonical representations}, and have finite image when $K_X+\Delta_X$ is semi-ample over $T$, see~\cite[Theorem 1.2, Corollary 1.3]{Fujino_RelativeLogPluricanonical}.

\begin{remark}\label{rmk: different invariants mu}
    Usually, we care about properties of $R_X=R(X/T,K_X+\Delta_X)$ that hold if and only if they are satisfied by some Veronese subalgebra, such as finite generation (see Lemma~\ref{lem: veronese finite generation}). Thus, the fact that $\Bim(X/T,\Delta_X)$ acts only on a Veronese subalgebra of $R_X$ does not represent an issue.
    However, there is a slight misalignment between the action defined in Lemma~\ref{lem: Bir acts on log canonical} and that in Lemma~\ref{lem: log canonical ring slc pair} in the case where $G=\langle\mu\rangle$ is given by the gluing involution associated to an slc pair. Indeed, in Lemma~\ref{lem: log canonical ring slc pair} the action is twisted by a factor of $(-1)^m$ in degree $m$. Hence the two actions differ only in odd degrees.
    We will thus always take $\ell$ to be even, which justifies the following definition. Note also that we use the different notation $R_X^{\mu}$ for the $\mu$-invariants under the ``twisted'' action, as $\mu$ is an element of $\langle\mu\rangle$ rather than the group itself.
\end{remark}

\begin{definition}\label{def: invariants under Bir}
    Let $(X,\Delta_X)$ be an lc pair, $f\colon X\rightarrow T$ a projective morphism to a complex analytic space $T$, and $R_X:=R(X/T,K_X+\Delta_X)$. Given a finite subgroup $G\subseteq\Bim(X/T,\Delta_X)$, we denote by $R_X^G$ the subalgebra of $R_X^{(\ell)}$ of $G$-invariants, where $\ell$ is the smallest positive \emph{even} integer such that $\ell(K_X+\Delta_X)$ is Cartier.
\end{definition}

Since we only care about finite generation of algebras, which is is invariant under passing to Veronese subalgebras by Lemma~\ref{lem: veronese finite generation}, we will not keep track of $\ell$ in the paper, in particular in Theorem~\ref{thm: equivariant sklt} and its proof.

\section{Good stable models}\label{sec: good stable models}

In this section we introduce a notion of \emph{stable model} for slc pairs, extending log canonical models to the non-normal case. We limit ourselves to treating the case where the exceptional locus of the map to the stable model does not contain any component of the conductor; for this reason, we call them \emph{good stable models}.
The general case is much more subtle, see for instance Example~\ref{exm: Proj RX not deminormal}.
We also provide examples of maps that should reasonably be considered stable models but fail to be good stable models, showing that the requirements~\eqref{property: defined up to codim 2} and~\eqref{property: birational conductor} in Definition~\ref{def: stable model slc pairs} are too restrictive in general.

\subsection{Definition and examples of good stable models}

A stable model is required to satisfy a large amount of properties, many of which will be proven to be redundant.

\begin{definition}[Good stable models]\label{def: stable model slc pairs}
    Let $(X,\Delta_X)$ be an slc pair, and let $(Y,\Delta_Y)$ be a demi-normal pair. Let $f\colon X\rightarrow T$, $g\colon Y\rightarrow T$ be projective morphisms, and let
    \[
    \begin{tikzcd}
        \phi\colon X\arrow[r,dashed] & Y
    \end{tikzcd}
    \]
    be a proper meromorphic map over $T$. We say that $\phi$ is a \emph{good stable model over $T$} if it satisfies the following properties:
    \begin{enumerate}
        \item\label{property: Y slc} $(Y,\Delta_Y)$ is slc;
        \item\label{property: ample} $K_Y+\Delta_Y$ is ample over $T$;
        \item\label{property: defined up to codim 2} the indeterminacy locus of $\phi$ contains no components of the conductor $D_X$ of $X$;
        \item\label{property: birational contraction} $\phi$ is a bimeromorphic contraction, that is, $\mathrm{codim}_Y(\Ex(\phi^{-1}))\geq2$;
        \item\label{property: normalization lc model} the induced map between normalization pairs
        \[
        \begin{tikzcd}
            \phibar\colon (\Xbar,D_{\Xbar}+\Delta_{\Xbar})\arrow[r,dashed] & (\Ybar,D_{\Ybar}+\Delta_{\Ybar})
        \end{tikzcd}
        \]
        is a log canonical model over $T$;
        \item\label{property: pushforward normalization} $\phibar_*D_{\Xbar}=D_{\Ybar}$, and $\phibar_*\Delta_{\Xbar}=\Delta_{\Ybar}$;
        \item\label{property: pushforward} $\phi_*\Delta_X=\Delta_Y$.
        \item\label{property: birational conductor} the restriction of $\phi$ to $D_X$ defines a bimeromorphic map $\phi|_{D_X}\colon D_X\dashrightarrow D_Y$. 
    \end{enumerate}
    When $X$ is normal, we call the stable model simply a log canonical model, as the two notions coincide.
\end{definition}

Note that the dependence on the base $T$ appears only in properties~\eqref{property: ample} and~\eqref{property: normalization lc model}. The proof of the following lemma is immediate, hence omitted.

\begin{lemma}\label{lem: redundant def stable model}
    Let $(X,\Delta_X)$ be an slc pair, and let $(Y,\Delta_Y)$ be a demi-normal pair with $\QQ$-Cartier log canonical divisor $K_Y+\Delta_Y$. Let $f\colon X\rightarrow T$, $g\colon Y\rightarrow T$ be projective morphisms, and let $\phi\colon X\dashrightarrow Y$ be a meromorphic map over $T$. Then, $\phi$ is a good stable model if and only if it satisfies properties \eqref{property: defined up to codim 2}, \eqref{property: normalization lc model}, and \eqref{property: birational conductor}.

    Moreover, if $\phi$ is a good stable model, then for every divisor $E$ over $X$, we have that
    \begin{equation}\label{eq: inequality discrepancies}
        a(E,X,\Delta_X)\leq a(E,Y,\Delta_Y).
    \end{equation}
    Conversely, if inequality \eqref{eq: inequality discrepancies} holds for the $\phi$-exceptional divisors and properties~\eqref{property: Y slc}, \eqref{property: ample}, \eqref{property: defined up to codim 2}, \eqref{property: birational contraction}, \eqref{property: pushforward} and \eqref{property: birational conductor} are satisfied, then $\phi$ is a good stable model.
\end{lemma}

Contrarily to what happens for normal spaces, the property of $\phi$ being regular on an open subset $U$ with complement of codimension at least 2 is not automatically satisfied, as $\phi$ could fail to be regular at the generic points of the conductor. Note that asking $D_X\dashrightarrow D_Y$ to be bimeromorphic implicitly requires $\phi$ to be defined at generic points of $D_X$.

\begin{example}\label{exm: bad conductor maps}
    Let us present some examples of morphisms $\phi$ that should reasonably be considered to be stable models even if they fail to be \emph{good} stable models. The examples will fail properties~\eqref{property: defined up to codim 2} or~\eqref{property: birational conductor}. Here, the base is $T=\Spec\CC$.
    \begin{enumerate}
        \item\label{not defined codim 1 example 1} We present a simple example where $\phi$ is not defined at the generic point of the conductor. Let $S$ be a smooth surface with $K_S$ ample. Let $p\in S$, and $\Bl_{p}S$ the blowup of $S$ along $p$, with exceptional $(-1)$-curve $E$. Let $X=\Bl_{p}S\cup_{E}\Bl_{p}S$ be the union of two copies of $S_1$ glued along $E$, thus $D_X=E$. Note that $X$ is simple normal crossing, and $K_X$ is big on $X$ as its pullback on each irreducible component $\Bl_{p}S$ is $K_{\Bl_{p}S}+E$. Note that $(K_{\Bl_{p}S}+E)|_E\simeq K_E\simeq\cO_{\PP^1}(-2)$, hence $R_X\simeq R_S\times R_S$. It follows that $R_X$ is finitely generated and $\phi\colon X\dashrightarrow \Proj R_X\simeq S\sqcup S$ satisfies all the properties of Definition~\ref{def: stable model slc pairs}, except for~\eqref{property: defined up to codim 2} and~\eqref{property: birational conductor}.
        \item\label{conductor example 1}
        Let $S$ and $p\in S$ be as in the previous example. Let $X$ be the variety obtained by gluing $\Bl_{p}S$ along an involution $\alpha$ of the exceptional divisor $E\simeq\PP^1$ with a unique fixed point.
        Note that $X$ is slc and Gorenstein, as it has only nodal singularities or pinch points along $D_X$. The induced birational morphism $\phi\colon X\rightarrow S$ automatically satisfies properties~\eqref{property: defined up to codim 2} and \eqref{property: pushforward} of Definition~\ref{def: stable model slc pairs}. Moreover, the normalization pair of $X$ is $(\Xbar=\Bl_{p}S,E)$, while $S$ is smooth, hence~\eqref{property: normalization lc model} is also satisfied. Therefore, $\phi$ has all the properties of a good stable model, except for~\eqref{property: birational conductor}.
        \item\label{conductor example 2} Let $Z$ be a smooth projective surface, and $C\subset Z$ a smooth curve with $K_Z+C$ ample. Let $Y$ be the union of two copies of $Z$ along $C$. Note that $Y$ is a stable slc pair with only SNC singularities and conductor $C$. Let $p\in C\subset Y$ and $X_0=\Bl_pY$. Then, $X_0=\Bl_pZ\cup_{C'}\Bl_pZ$, where $C'$ is the strict transform of $C$; note that $C'\simeq C$. Let $X_1$ be the scheme obtained by gluing $X_0$ along the two exceptional divisors $E\subset Z$, whose image in $X_1$ is still denoted by $E$. Then, $(X_1,0)$ is slc with conductor $E\cup C'$. Let $\phi_i\colon X_i\rightarrow Y$ be the induced morphisms. Then, a simple computation shows that both $\phi_0$ and $\phi_1$ satisfy all properties of Definition~\ref{def: stable model slc pairs} except that $\phi_1$ fails~\eqref{property: birational conductor}. Note that $\phi_1$ contracts the component $E$ of the conductor, even though it restricts surjectively to $D_{X_1}=E\cup C'\rightarrow C=D_Y$.
    \end{enumerate}
\end{example}

\begin{remark}\label{rmk: non stable models}
    Let $(X,\Delta_X)$ be a non-normal slc pair with ample log canonical divisor. Then, the normalization map $\pi_X\colon\Xbar\longrightarrow X$ is not a stable model for any boundary on $\Xbar$, as it fails properties \eqref{property: birational contraction} and \eqref{property: pushforward normalization} of Definition~\ref{def: stable model slc pairs}. Similarly, the inverse $\pi_X^{-1}\colon X\dashrightarrow\Xbar$ of the normalization map fails properties \eqref{property: birational contraction} and \eqref{property: pushforward normalization} of Definition~\ref{def: stable model slc pairs}.
\end{remark}

\begin{remark}\label{rmk: general definition stable model}
    Example~\ref{exm: bad conductor maps} shows that to obtain a possible general notion of stable model we need to remove condition~\eqref{property: defined up to codim 2}, and of course~\eqref{property: birational conductor}. Remark~\ref{rmk: non stable models} suggests that all other properties are to be maintained.
    We omit a further discussion on this, as it is beyond the scope of this work.
\end{remark}

\subsection{Good stable models as Projan of log canonical algebras}\label{subsec: good stable models as Projan}

We prove that good stable models are unique, and coincide with the Projan of the log canonical algebra.

\begin{proposition}\label{prop: Y=proj RX good stable model}
    Let $f\colon X\rightarrow T$ be a projective morphism between complex analytic spaces, $(X,\Delta_X)$ an slc pair, and $\phi\colon X\dashrightarrow Y$ a good stable model over $T$. Then, $\phibar|_{D_{\Xbar}^n}\colon D_{\Xbar}^n\dashrightarrow D_{\Ybar}^n$ induces a natural injection $R_{D_{\Ybar}^n}\hookrightarrow R_{D_{\Xbar}^n}$ between log canonical algebras, that is equivariant with respect to the gluing involutions and makes
    \begin{equation}\label{diag: commutative restriction}
    \begin{tikzcd}
        R_{D_{\Xbar}^n}\arrow[r,hookleftarrow] & R_{D_{\Ybar}^n}\\
        R_{\Xbar}\arrow[u,"\res_n"] & R_{\Ybar}\arrow[l,"\simeq"]\arrow[u,"\res_n"]
    \end{tikzcd}
    \end{equation}
    commutative.
    In particular, $R_Y\simeq R_X$ and the good stable model $Y\simeq\Projan_T R_X$ is unique.
\end{proposition}
\begin{proof}
    The statement is local on $T$, hence we are allowed to restrict $T$ as needed. In particular, we can assume $T$ to be Stein and $\phibar:\Xbar\dashrightarrow\Ybar$ to admit a common resolution $\Xbar\xleftarrow{p}V\xrightarrow{q}\Ybar$, eventually after restricting to a relatively compact open subset of $T$. Write
    \[
        E=p^*(K_{\Xbar}+D_{\Xbar}+\Delta_{\Xbar})-q^*(K_{\Ybar}+D_{\Ybar}+\Delta_{\Ybar}),
    \]
    which is effective and $q$-exceptional by the properties of log canonical models. Let $Z\subset V$ be the strict transform of $D_{\Xbar}$ along $p$, which coincides with the strict transform of $D_{\Ybar}$ along $q$ by Definition~\ref{def: stable model slc pairs}\eqref{property: birational conductor}. Let $Z^n$ be the normalization of $Z$. We obtain the following commutative diagram:
    \begin{equation}
    \begin{tikzcd}
        & Z^n\arrow[dd]\arrow[ld,"p_D"']\arrow[rd,"q_D"]\\
        D_{\Xbar}^n\arrow[dd]\arrow[rr,dashed, crossing over] & & D_{\Ybar}^n\arrow[dd]\\
        & V\arrow[ld,"p"']\arrow[rd,"q"]\\
        \Xbar\arrow[rr,dashed] & & \Ybar
    \end{tikzcd}
    \end{equation}
    Since $K_{\Xbar}+D_{\Xbar}+\Delta_{\Xbar}$ and $K_{\Ybar}+D_{\Ybar}+\Delta_{\Ybar}$ are $\mathbb{Q}$-Cartier, pulling them back to $V$ and restricting to $Z^n$ yields the same result as taking the different first and then pulling back to $Z^n$. Restricting the divisor equality for $E$ to $Z^n$ thus yields:
    \[
        p_D^*(K_{D_{\Xbar}^n}+\Diff_{D_{\Xbar}^n}(\Delta_{\Xbar})) = q_D^*(K_{D_{\Ybar}^n}+\Diff_{D_{\Ybar}^n}(\Delta_{\Ybar})) + E_Z,
    \]
    where $E_Z := E|_{Z^n}$ remains effective but is no longer $q_D$-exceptional in general. 
    
    Since $\phi$ is a good stable model, the map between conductors is bimeromorphic, hence $p_D$ and $q_D$ are bimeromorphic morphisms between normal complex analytic spaces. Let $f_D:D_{\Xbar}^n\rightarrow T$ and $g_D:D_{\Ybar}^n\rightarrow T$ be the projective morphisms induced by the structure morphisms $f:X\rightarrow T$ and $g:Y\rightarrow T$, respectively. Since $E_Z\geq 0$, pushing forward to $T$ for all sufficiently divisible powers $m\geq 0$ yields a natural inclusion
    \begin{align*}
        g_{D*}\cO_{D_{\Ybar}^n}(m(K_{D_{\Ybar}^n}+\Diff_{D_{\Ybar}^n}(\Delta_{\Ybar}))) & \simeq g_{D*}q_{D*}\cO_{Z^n}(m q_D^*(K_{D_{\Ybar}^n}+\Diff_{D_{\Ybar}^n}(\Delta_{\Ybar})))\\
        & \subseteq g_{D*}q_{D*}\cO_{Z^n}(m p_D^*(K_{D_{\Xbar}^n}+\Diff_{D_{\Xbar}^n}(\Delta_{\Xbar})))\\
        & \simeq f_{D*}p_{D*}\cO_{Z^n}(m p_D^*(K_{D_{\Xbar}^n}+\Diff_{D_{\Xbar}^n}(\Delta_{\Xbar})))\\
        & \simeq f_{D*}\cO_{D_{\Xbar}^n}(m(K_{D_{\Xbar}^n}+\Diff_{D_{\Xbar}^n}(\Delta_{\Xbar}))).
    \end{align*}
    Summing over all $m$ and applying~\cite[Lemma 2.19]{Kollar-singularities}, we obtain the desired injection of graded rings $R_{D_{\Ybar}^n}\hookrightarrow R_{D_{\Xbar}^n}$, which makes diagram~\eqref{diag: commutative restriction} commute. The inclusion is $\mu$-equivariant, since the bimeromorphic map $D_{\Xbar}^n\dashrightarrow D_{\Ybar}^n$ is $\mu$-equivariant, being induced by $\phi\colon X\dashrightarrow Y$.

    The global isomorphism $R_X\simeq R_Y$ follows from diagram~\eqref{diag: commutative restriction} and Lemma~\ref{lem: log canonical ring slc pair}.
    Since $K_Y+\Delta_Y$ is ample by Definition~\ref{def: stable model slc pairs}\eqref{property: ample}, we conclude that $Y\simeq\Proj R_Y\simeq\Proj R_X$.
\end{proof}

Let $(X,\Delta_X)$ be an slc pair projective over a complex analytic space $T$, suppose the log canonical ring $R_X$ to be locally finitely generated, and that $X\dashrightarrow Y:=\Projan_T R_X$ is bimeromorphic. Let $\Delta_Y=\phi_*\Delta_X$. If $X$ is normal, then it is well-known that $\Proj R_X$ is the log canonical model of $X$; see for instance \cite{Reid-canonical} or~\cite[Theorems 1.15, 1.26]{Kollar-singularities} for the algebraic case.
However, for an slc pair it is not even clear that $Y:=\Proj R_X$ is demi-normal, which indeed fails by Example~\ref{exm: Proj RX not deminormal} below. Even then, the relationship between the finite generation of $R_X$ and $R_{\Xbar}$ is somewhat obscure, as the extension $R_X\hookrightarrow R_{\Xbar}$ is not always finite. Indeed, there are cases where $R_{\Xbar}$ is a finitely generated ring but $R_X$ is not, even when $X$ is a normal-crossing surface, see~\cite{Kollar_examples_SNC}.

We prove that $\Projan_T R_X$ is a good stable model if the map $X\dashrightarrow\Projan_T R_X$ is holomorphic in codimension 1 and does not contract any component of the conductor of $X$.

\begin{theorem}\label{thm: projRX is stable model}
    Let $(X,\Delta_X)$ be an slc pair projective over a semi-normal complex analytic space $T$, such that each irreducible component $X_i$ of $X$ dominates some irreducible component $T_i$ of $T$ (not necessarily distinct). Denote by $\overline{X}_i$ the normalization of $X_i$ and by $(\overline{X}_i,D_{\overline{X}_i}+\Delta_{\overline{X}_i})$ the induced pair. Suppose that the $\cO_T$-algebra $R_X:=R(X/T,K_X+\Delta_X)$ is locally finitely generated, and that $K_{\overline{X}_i}+D_{\overline{X}_i}+\Delta_{\overline{X}_i}$ on $\overline{X}_i$ is relatively big over $T_i$ for all $i$. Then the induced map $\phi\colon X\dashrightarrow\Projan_T R_X$ is bimeromorphic.
    
    Suppose further that $\phi$ is holomorphic in codimension 1 and does not contract any component of the conductor $D_X$. Then, $\phi$ is the unique good stable model of $X$ over $T$, and it induces an identification of log canonical algebras.
    Moreover, if $K_X+\Delta_X$ is relatively nef over $T$, then $\phi$ is a morphism.
    
    In particular, $R_{\Xbar}:=R(\Xbar/T,K_{\Xbar}+D_{\Xbar}+\Delta_{\Xbar})$ is also finitely generated and the induced map $\Projan_T R_{\Xbar}\dashrightarrow\Projan_T R_X$ is a morphism that coincides with the normalization of $\Projan_T R_X$.
\end{theorem}

First, we prove a proposition in a more general setting.

\begin{proposition}\label{prop: properties of Proj deminormal}
    Let $X$ be a demi-normal complex analytic space, with a projective morphism $f\colon X\rightarrow T$ to a semi-normal complex analytic space $T$, such that each irreducible component $X_i$ of $X$ dominates some irreducible component $T_i$ of $T$ (not necessarily distinct). Let $L$ be an effective Cartier divisor on $X$.
    Suppose that $R(X/T,L)$ is generated by $f_*\cO_X(L)$ as an $\cO_T$-algebra, and that $L_i:=L|_{\overline{X}_i}$ on the normalization $\overline{X}_i$ of $X_i$ is relatively big over $T_i$ for all $i$.
    Let $\phi\colon X\dashrightarrow\PP_T^N$ be the map given by $|L|_T$, defined after restricting $T$ to any relatively compact Stein open subset. Let $Y$ be the closure of the image and $|H|_T$ the hyperplane class on $Y$. Suppose that $\Ex(\phi)$ does not contain any component of the conductor $D_X$. Then:
    \begin{enumerate}
        \item $\phi$ is holomorphic in codimension 1 and $\phi$ is bimeromorphic;
        \item $Y$ is demi-normal;
        \item $\phi_*^{-1}D_Y=D_X$ and $\phi_*D_X=D_Y$;
        \item $Y\setmin\phi(X\setmin\Bs|L|_T)$ has codimension $\geq2$;
        \item every divisor in $\Bs|L|_T$ is contracted by $\phi$;
        \item $\phi_*|L|_T=|H|_T$;
        \item if $L$ is nef over $T$, then $|L|_T$ is base-point free over $T$.
    \end{enumerate}
\end{proposition}
\begin{proof}
    The proof is local on $T$, therefore we can replace $T$ with a relatively compact Stein open subset. We denote by $g:Y\rightarrow T$ the projective structure map.
    
    By assumption, $\phi$ is holomorphic on a Zariski-open $U$ of $X$ whose complement has codimension at least 2. Moreover, the relative bigness of $L_i$ on $\overline{X}_i$ immediately implies that $\phi$ is bimeromorphic.

    Now, let $X\xleftarrow{p}V''\xrightarrow{q}Y$ be the closure of the graph of $\phi$, with the two projections. Since $\phi$ is defined on the open $U$ whose complement has codimension at least 2, $p\colon V''\rightarrow X$ is an isomorphism over $U$ and $p^{-1}(U)\subset V''$ is at most nodal in codimension 1. However, it is possible for $V''$ to have other singularities in codimension 1 that get contracted to a subscheme of $X\setminus U$. To fix this, we first take a semi-normalization of $V''$ and then blowup the singular codimension 1 loci outside $p^{-1}(U)$. The resulting complex analytic space $V'$ is now nodal in codimension 1, thus we can take its demi-normal modification $V\rightarrow V'$ (Definition~\ref{def: deminormal modification}). We still use $p$ and $q$ to denote the maps from $V$ to $X$ and $Y$, respectively. Then, $p\colon V\rightarrow X$ is a bimeromorphic morphism between demi-normal spaces, induces an isomorphism $p^{-1}(U)\rightarrow U$, and all singular points of codimension 1 in $X$ (resp. $V$) are contained in $U$ (resp. $p^{-1}(U)$). It follows that $p_*\cO_V\simeq\cO_X$. Therefore, if we set $M:=p^*L$, the pair $(V,M)$ satisfies the same assumptions as in the statement, and $\Bs|M|_T=p^{-1}\Bs|L|_T$. It is then enough to prove the statement assuming that $\phi$ is a surjective proper morphism, and
    \begin{equation}\label{eq: relation linear systems}
        |L|_T=\phi^*|H|_T+F,
    \end{equation}
    where $F$ is a Cartier divisor with support equal to the base locus of $|L|_T$. In particular, $g_*\cO_Y(mH)\simeq f_*\cO_{X}(mL)$ for all $m\geq0$.

    Now, let $X\xrightarrow{\phi'}Y'\xrightarrow{\psi}Y$ be any factorization of $\phi$ with $\phi'$ a birational morphism and $\psi$ finite. Then, $m\psi^*H=\psi^*(mH)$ is very ample over $T$ for $m\gg0$, and
    \[
    g_*\cO_{Y}(mH)\subseteq g_*\psi_*\cO_{Y'}(m\psi^*H)\subseteq f_*\cO_{X}(mL),
    \]
    as the morphisms are dominant and satisfy an equality between linear systems as in~\eqref{eq: relation linear systems}. Since the composite of inclusions is an equality, we then get $Y'\simeq Y$. It follows that any such factorization is the trivial one, in particular $\phi_*\cO_X=\cO_Y$ and the fibers are connected.
    
    We are ready to show that $Y$ is demi-normal. Since the fibers of $\phi$ are connected, there is a Zariski-open subset of $Y$ over which $\phi$ is an isomorphism, and whose complements has codimension at least 2. In particular, $Y$ is nodal in codimension 1. Then, Example~\ref{exm: deminormal modifications}\eqref{exm: deminormal stein factorization} shows that $Y$ is demi-normal.
    
    Note that all Weil divisors $B$ in $X$ not contained in $D_X$ are either contracted or mapped to a Weil divisor in $Y$ not contained in $D_Y$, and that all irreducible components of $D_X$ are either contracted or mapped to an irreducible component of $D_Y$.
    Since the components of the conductor are not contracted by assumption, we get that $\phi_*^{-1}D_Y=D_X$ and $\phi_*D_X=D_Y$.
    
    The same argument as in~\cite[Proposition 1.16]{Kollar-singularities} shows that if a Weil divisor $B$ not contained in $D_X$ is not contracted then $B\not\subset F=\Bs|L|_T$. Since no component of $D_X$ is contained in the base locus, it follows that $\phi(\Bs|L|_T)$ has codimension at least 2 in $Y$, and $\phi_*|L|_T=|H|_T$.

    The fact that $|L|$ is base-point free when $L$ is nef follows from the next lemma.
\end{proof}

\begin{lemma}
    Let $\phi\colon X\rightarrow Y$ be a projective bimeromorphic morphism between demi-normal complex analytic spaces, and $F$ an effective $\QQ$-Cartier $\QQ$-divisor that is $\phi$-nef and contracted by $\phi$. Then, $F=0$.
\end{lemma}
\begin{proof}
    The result is standard when $X$ and $Y$ are normal, and it follows from the Negativity Lemma (\cite[Lemma 3.39]{Kollar_Mori_1998}, \cite[Corollary 2.16]{Enokizono-Hashizume-analytic-MMP}). In general, let $\phibar\colon \Xbar\rightarrow\Ybar$ be the induced morphism between the normalizations, and note that $\pi_X^*F$ is again effective, $\phibar$-nef, and contracted by $\phibar$.
    Therefore, $\pi_X^*F=0$, thus $F=0$.
\end{proof}

\begin{proof}[Proof of Theorem~\ref{thm: projRX is stable model}]
    The proof is local on $T$, therefore we can replace $T$ with a relatively compact Stein open subset, further restricting it where needed, and assume $R_X$ to be finitely generated.
    
    Since $R_X$ is finitely generated, there is an integer $\ell>0$ such that $\ell(K_X+\Delta_X)$ is Cartier and the $\ell$-th Veronese subalgebra of $R_X$ is generated by $f_*\cO_{X}(\ell(K_X+\Delta_X))$. Then, $L:=\ell(K_X+\Delta)$ satisfies the assumptions of Proposition~\ref{prop: properties of Proj deminormal}.
    Therefore, there is a natural bimeromorphic map $\phi\colon X\dashrightarrow Y$ which is holomorphic in codimension 1 and with demi-normal target, together with a very ample Cartier divisor $H$ on $Y$ satisfying $H\sim\phi_*L$, and $f_*\cO_{X}(mL)\simeq g_*\cO_{Y}(mH)$ via $\phi$, for all $m\geq0$. Moreover, $\Delta_Y:=\phi_*\Delta_X$ is a Weil divisor with no components contained in $D_X$, since $\phi^{-1}_*D_Y=D_X$. Then,
    \[
        \ell(K_Y+\Delta_Y)\sim\phi_*\ell(K_X+\Delta_X)=\phi_*L\sim H,
    \]
    hence $K_Y+\Delta_Y$ is $\QQ$-Cartier and ample over $T$. Since $\phi$ does not contract the conductor by assumption, by Lemma~\ref{lem: redundant def stable model} we are only left with showing that $\phibar$ is a log canonical model.

    For this, we can reduce to the case where $\phi$ is a morphism, by taking the demi-normal modification of a partial resolution of the closure of the graph of $\phi$, as done in the proof of Proposition~\ref{prop: properties of Proj deminormal}. Recall that we have
    \[
        \phi^*(\ell(K_Y+\Delta_Y))+F=\ell(K_X+\Delta_X)
    \]
    with $F$ effective and $\phi$-exceptional. Pulling back to the normalization, we get
    \[
        \phibar^*(\ell(K_{\Ybar}+D_{\Ybar}+\Delta_{\Ybar}))+\pi_X^*F=\ell(K_{\Xbar}+D_{\Xbar}+\Delta_{\Xbar}).
    \]
    Then, $\overline{g}_*\cO_{\Ybar}(m(K_{\Ybar}+D_{\Ybar}+\Delta_{\Ybar}))=\overline{f}_*\cO_{\Xbar}(m(K_{\Xbar}+D_{\Xbar}+\Delta_{\Xbar}))$ for all $m\geq0$ sufficiently divisible. It follows that $R_{\Xbar}=R(\Xbar/T,K_{\Xbar}+D_{\Xbar}+\Delta_{\Xbar})$ is finitely generated and $\Proj R_{\Xbar}=\Ybar$ is indeed the log canonical model. Moreover, by Proposition~\ref{prop: properties of Proj deminormal} we know that $\phibar_*D_{\Xbar}=D_{\Ybar}$, while $\phi_*\Delta_{\Xbar}=\Delta_{\Ybar}$ holds by definition of $\Delta_{Y}$.
\end{proof}

\begin{example}\label{exm: Proj RX not deminormal}
    We provide an example showing that $\Proj R_X$ can fail to be demi-normal if the exceptional locus of the natural map $X\dashrightarrow \Proj R_X$ contains a component of $D_X$.

    Let $S$ be a minimal smooth projective surface containing a smooth elliptic curve $E$. Let $T$ be the log canonical model of $(S,E)$, which is obtained by contracting $E$ to a simple elliptic singularity $p\in T$. Let $X$ be the union of two copies of $S$ glued along $E$, and let $R_X:=R(X,K_X)$. Note that the restriction map
    \[
    \begin{tikzcd}
        \res\colon R_S:=R(S,K_S+E)\arrow[r] & R(E,K_E)=:R_E\simeq\CC[t]
    \end{tikzcd}
    \]
    is surjective by~\cite[Proposition 4.45]{Kollar_Mori_1998} or Kawamata-Viehweg vanishing. Moreover, for every $m\geq0$ there exists an exact sequence
    \[
    \begin{tikzcd}[column sep=small]
        0\arrow[r] & \H^0(X,mK_X)\arrow[r] & \H^0(S,m(K_S+E))\oplus \H^0(S,m(K_S+E))\arrow[r,"\phi"] & \H^0(E,mK_E)\arrow[r] & 0
    \end{tikzcd}
    \]
    where $\phi(s_1,s_2)=\res(s_1)-(-1)^m\res(s_2)$. In particular, $R_X\simeq R_S\times_{\CC[t]}R_S$.
    From this and the surjectivity of $\res$, it follows that $R_X$ is finitely generated and $\Proj R_X$ is the union of two copies of $T$ glued along the elliptic singularity $p$, hence not demi-normal.
\end{example}

An important question is whether a good stable model of a semi-dlt pair induces another good stable model for the conductor.
This fails in general, but there are useful situations where the answer is positive.

\begin{lemma}\label{lem: if birational conductor then stable model}
    Let $\phi\colon (X,\Delta_X)\rightarrow(Y,\Delta_Y)$ be a good stable model over $T$ of a semi-dlt pair $(X,\Delta_X)$ projective over a complex analytic space $T$. Suppose further that $\phi$ is a crepant bimeromorphic morphism, that is, $K_X+\Delta_X=\phi^*(K_Y+\Delta_Y)$, and that $\Ex(\phi)$ does not contain any stratum of the conductor $D_X$. Let $\phi_D\colon D_{\Xbar}\rightarrow\widetilde{D}_{\Ybar}$ be the morphism induced by taking the Stein factorization of $D_{\Xbar}\rightarrow D_{\Ybar}$. Then, $\widetilde{D}_{\Ybar}$ is demi-normal, $\phi_{D_{\Xbar}}$ is a crepant good stable model over $T$, and $\Ex(\phi_D)$ does not contain any stratum of the conductor. Moreover, $\widetilde{D}_{\Ybar}\simeq\Projan_T R_{D_{\Xbar}}$.
\end{lemma}

This is a special case of the following more general lemma.

\begin{lemma}\label{lem: birational coeff 1 boundary then stable model}
    Let $(X,\Delta_X)$ be a dlt pair that is projective over a complex analytic space $T$. Let $\phi\colon (X,\Delta_X)\rightarrow(Y,\Delta_Y)$ be a log canonical model over $T$ that is also a crepant bimeromorphic morphism, meaning $K_X+\Delta_X=\phi^*(K_Y+\Delta_Y)$. Let $S_X\leq\lfloor\Delta_X\rfloor$ be a reduced subdivisor of the coefficient 1 part of the boundary, such that $\phi|_{S_X}\colon S_X\rightarrow Y$ is bimeromorphic onto its image $S_Y\leq\lfloor\Delta_Y\rfloor$ and does not contract any stratum of $S_X$. Let $\phi_1\colon S_X\rightarrow\widetilde{S}_{Y}$ be the morphism induced by taking the Stein factorization of $S_X\rightarrow S_Y$. Then, $(S_X,\Diff_{S_X}(\Delta_X-S_X))$ is semi-dlt, $\widetilde{S}_{Y}$ is demi-normal, $\phi_1$ is a crepant good stable model over $T$, and $\Ex(\phi_1)$ does not contain any stratum of $S_X$. Moreover, $\widetilde{S}_{Y}\simeq\Projan_T R_{S_X}$.
\end{lemma}
\begin{proof}
    By the standard (semi-)dlt adjunction in Lemma~\ref{lem: semidlt adjunction} and~\cite[Corollary 4.5]{fujino2016vanishing}, the restriction to the subdivisor $S_X$ yields a well-defined semi-dlt pair $(S_X, \Diff_{S_X}(\Delta_X-S_X))$. In particular, $S_X$ is a demi-normal complex analytic space.
    
    Taking the Stein factorization of the proper morphism $S_X\rightarrow S_Y$ yields a morphism $\phi_1\colon  S_X \rightarrow\widetilde{S}_Y$ with connected fibers. Then, $\widetilde{S}_Y$ is demi-normal because the domain $S_X$ is demi-normal, the map does not contract any component of the conductor of $S_X$, and $S_Y$ is nodal in codimension 1, see Example~\ref{exm: deminormal modifications}\eqref{exm: deminormal stein factorization}. Restricting the crepant identity $K_X+S_X+(\Delta_X-S_X)=\phi^*(K_Y+\Delta_Y)$ to $S_X$ via adjunction (or \cite[Proposition 4.6]{Kollar-singularities}) shows that $\phi_1$ is a crepant morphism of pairs.
    Since $\phi_1$ is a crepant bimeromorphic morphism with connected fibers from a semi-dlt pair, the target pair $(\widetilde{S}_Y, \Delta_{\widetilde{S}_Y})$ is automatically slc. Furthermore, since $K_Y+\Delta_Y$ is ample over $T$ by the definition of a log canonical model, and $\widetilde{S}_Y$ is finite over $Y$, $(\widetilde{S}_Y,\Delta_{\widetilde{S}_Y})$ is a stable slc pair.
    Since $\phi_1$ is a crepant bimeromorphic morphism from a semi-dlt pair to a target polarized by an ample log canonical class, it satisfies every condition of a good stable model over $T$. The isomorphism $\widetilde{S}_{Y}\simeq\Proj R_{S_X}$ follows from Proposition~\ref{prop: Y=proj RX good stable model}.
\end{proof}

\section{Finite generation of the log canonical algebra}\label{sec: finite generation}

In this section, we prove Theorems~\ref{thm: general version} and~\ref{thm: equivariant sklt}, starting from the latter.
First, we note that we can perform the following reductions.

\begin{remark}[Reductions]\label{rmk: reductions}
    Since the statement of Theorem~\ref{thm: equivariant sklt} is local on $T$, we can assume that we have restricted $T$ to a Stein open subset that is relatively compact. In particular, we are allowed to apply any result from resolution of singularities that we need. Moreover, we can assume that every relevant coherent sheaf on $T$ is finitely globally generated, and that $X$ carries an $f$-very ample divisor $A_X$.
    In the proofs, we will also need to apply some results from the minimal model program for projective morphisms between complex analytic spaces, whose existence has only been established around Stein compact subsets satisfying property $\cP$ (\cite{Fujino_MMPanalytic,Enokizono-Hashizume-analytic-MMP}). On the other hand, since we are free to arbitrarily restrict $T$, we can always assume that all required results from MMP work globally over $T$. To simplify the exposition, we will often tacitly perform all these reductions and avoid naming the compact subset of $T$ we are working around, occasionally referring to this remark when needed.
    Finally, by Lemma~\ref{lem: veronese finite generation}, to prove the (local) finite generation of relative log canonical algebras we are allowed to take arbitrary $\ell$-Veronese subalgebras. Again, we will often do so without mentioning it, and take $\ell$ to be even because of Remark~\ref{rmk: different invariants mu}.
\end{remark}

\subsection{Proof of Theorem~\ref{thm: equivariant sklt}}\label{subsec: proof equivariant sklt}

We start by recalling the setting. Let $(X,\Delta_X)$ be a semi-klt pair, $T$ a semi-normal complex analytic space, and $f\colon X\rightarrow T$ a surjective projective morphism. Let $T^0\subset T$ be a Zariski-open subspace, meaning that its complement is an analytical subset, and set $X^0=T^0\times_{T}X$, $\Delta_{X^0}=\Delta_X|_{X^0}$. We may also assume that $X$ is not normal, hence $D_X\not=\emptyset$. We assume that $X^0$ intersects all log centers of $(X,\Delta_X)$ lying on $D_X$; in particular, $X^0$ intersects all strata of the conductor $D_X$ and all irreducible components of $X$.
Suppose that $(X^0,\Delta_{X^0})$ admits a stable model
\[
\begin{tikzcd}
    \phi^{0,c}\colon (X^0,\Delta_{X^0})\arrow[r,dashed] & (X^{0,c},\Delta_{X^{0,c}})
\end{tikzcd}
\]
whose exceptional locus does not contain any stratum of $D_{X^0}$, and that the normalization pair $(\Xbar,D_{\Xbar}+\Delta_{\Xbar})$ admits a log canonical model $\phibar^c\colon \Xbar\rightarrow(\Xbar^c,D_{\Xbar^c}+\Delta_{\Xbar^c})$. Let $G\subset\Bim(\Xbar/T,D_{\Xbar}+\Delta_{\Xbar})$ be a finite subgroup. Let us set $R_X:=R(X/T,K_X+\Delta_X)$, and similarly for $R_{\Xbar}$.
Our goal is then to show that the extension of $\cO_T$-algebras
\[
\begin{tikzcd}
    R_X\cap R_{\Xbar}^G\arrow[r,hookrightarrow] & R_{\Xbar}
\end{tikzcd}
\]
is finite, with the convention of Definition~\ref{def: invariants under Bir}.

Note that in the proof we are allowed to perform the reductions explained in Remark~\ref{rmk: reductions}; in particular, we will assume that $T$ is a sufficiently small Stein space, and we will tacitly pass to Veronese subalgebras of even degree when needed.


The proof is by induction on the dimension of $X$. The base case is for 0-dimensional pairs, thus trivially true. For the inductive step, we need the following lemma, that highlights the nice properties specific to semi-klt pairs.

\begin{lemma}\label{lem: conductor for sklt induction}
    With the above assumptions, and after suitably restricting $T$ as explained in Remark~\ref{rmk: reductions}, there exists a dlt pair $(\Xbar^m,D_{\Xbar^m}+\Delta_{\Xbar^m})$
    fitting in a commutative diagram of meromorphic maps over $T$
    \[
    \begin{tikzcd}
        D_{\Xbar}^n\arrow[d]\arrow[r,dashed,"\phibar_{D^n}^m"] & D_{\Xbar^m}^n\arrow[d]\arrow[r,"\phibar_{D^n}"] & D_{\Xbar^c}^n\arrow[d]\\
        D_{\Xbar}\arrow[r,dashed,"\phibar_D^m"]\arrow[d,hookrightarrow] & D_{\Xbar^m}\arrow[r,"\phibar_D"]\arrow[d,hookrightarrow] & D_{\Xbar^c}\arrow[d,hookrightarrow]\\
        \Xbar\arrow[d,"\pi_X"]\arrow[r,dashed,"\phibar^m"] & \Xbar^m\arrow[r,"\phibar"] & \Xbar^c\\
        X
    \end{tikzcd}
    \]
    where solid arrows are morphisms, such that the following hold.
    \begin{enumerate}[itemsep=1em]
        \item\label{point 1} $\phibar^m$ is the minimal model of $(\Xbar,D_{\Xbar}+\Delta_{\Xbar})$, and $\phibar^c=\phibar\circ\phibar^m$ the log canonical model, with $\phibar_*^m D_{\Xbar}=D_{\Xbar^m}$ and $\phibar_* D_{\Xbar^m}=D_{\Xbar^c}$.
        \item\label{point 2} $\phibar_D$, $\phibar_D^m$ and $\phibar_D^c:=\phibar_D\circ\phibar_D^m$ are bimeromorphic, and their exceptional loci do not contain any stratum of their conductors. The same holds for their normalizations $\phibar_{D^n}$, $\phibar_{D^n}^m$ and $\phibar_{D^n}^c:=\phibar_{D^n}\circ\phibar_{D^n}^m$.
        \item\label{point 3} $\lfloor\Delta_{\Xbar^m}\rfloor=0$ and $\lfloor\Diff_{D_{\Xbar^m}}(\Delta_{\Xbar^m})\rfloor=0$.
        \item\label{point 4} $\phibar_D$ is crepant and a good stable model of the semi-klt pair $(D_{\Xbar^m},\Diff_{D_{\Xbar^m}}(\Delta_{\Xbar^m}))$. The same holds for $\phibar_{D^n}$. In particular, $D_{\Xbar^c}$ is demi-normal, and $(D_{\Xbar^c},\Diff_{D_{\Xbar^c}}(\Delta_{\Xbar^c}))$ is a stable slc pair.
        \item\label{point 5} Let $R_{D_{\Xbar}^n}:=R(D_{\Xbar}^n/T,K_{D_{\Xbar}^n}+\Diff_{D_{\Xbar}^n}(\Delta_{\Xbar}))$, and similarly for $R_{D_{\Xbar^m}^n}$ and $R_{D_{\Xbar^c}^n}$. Then,
        \[
            R_{D_{\Xbar}^n}\supset R_{D_{\Xbar^m}^n}\simeq R_{D_{\Xbar^c}^n}.
        \]
    \end{enumerate}
\end{lemma}
\begin{proof}
    The existence of the diagram and point~\eqref{point 1} follow from Lemma~\ref{lem: existence minimal model if log canonical}.

    Point~\eqref{point 2} follows from the assumption that $X^0$ intersects all strata of the conductor of $X$, which are not contracted by $\phi^{0,c}$. Point~\eqref{point 3} follows from the equality $\phibar_*^m\Delta_{\Xbar}=\Delta_{\Xbar^m}$ and Lemma~\ref{lem: semidlt adjunction}.

    To prove~\eqref{point 4}, it is enough to show that $\phibar_D$ is a good stable model, as $\phibar_D$ is crepant (see for instance~\cite[Proposition 4.6]{Kollar-singularities} or Lemma~\ref{lem: birational coeff 1 boundary then stable model}). By Lemma~\ref{lem: if birational conductor then stable model}, $D_{\Xbar^m}$ has finitely generated log canonical algebra $R_{D_{\Xbar^m}}$, whose Proj is a stable model $Z$ that is finite over $D_{\Xbar^c}$. By Lemma~\ref{lem: surjectivity restriction minimal model}, the restriction map $R_{\Xbar^m}^{(\ell)}\twoheadrightarrow R_{D_{\Xbar^m}}^{(\ell)}$ is surjective for $\ell$ sufficiently divisible. It follows that $\Proj R_{D_{\Xbar^m}}\rightarrow\Proj R_{\Xbar^m}=\Xbar^c$ is a closed immersion. In particular, $Z=\Proj R_{D_{\Xbar^m}}\rightarrow D_{\Xbar^c}$ is an isomorphism, being a morphism of reduced subspaces with the same underlying set.

    The isomorphism in~\eqref{point 5} follows from $\phibar_{D^n}$ being a log canonical model, while the inclusion on the left follows from Proposition~\ref{prop: Y=proj RX good stable model}.
\end{proof}

We are ready to prove the inductive step for Theorem~\ref{thm: equivariant sklt}; we use the notation of Lemma~\ref{lem: conductor for sklt induction}.

\begin{proof}[Proof of Theorem~\ref{thm: equivariant sklt}]
Throughout, we assume $T$ to be Stein, see Remark~\ref{rmk: reductions}.

First, note that there is an honest $G$-action on the canonical model $(\Xbar^c,D_{\Xbar^c}+\Delta_{\Xbar^c})$, as $\Xbar^c=\Projan_T R_{\Xbar}$ and $G$ acts on a Veronese subalgebra of $R_{\Xbar}$ by Lemma~\ref{lem: Bir acts on log canonical}.
This action restricts to a $G$-action on $(D_{\Xbar^c}^n,\Diff_{D_{\Xbar^c}^n}(\Delta_{\Xbar^c}))$. Moreover, there is a canonical involution $\mu$ on this pair, defined as follows. Recall that $(\Xbar^{0,c},D_{\Xbar^{0,c}}+\Delta_{\Xbar^{0,c}})$ is the normalization of the stable model $(X^{0,c},\Delta_{X^{0,c}})$ of $X$ over $T^0$. Therefore, there is a gluing involution $\mu^0$ on $(D_{\Xbar^c}^n,\Diff_{D_{\Xbar^c}^n}(\Delta_{\Xbar^c}))|_{T^0}$, see \S\ref{subsec: normal and demi-normal modification}. Since all log centers on $D_{\Xbar^c}$ intersect $\Xbar^{0,c}$ by assumption (and~\cite[Proposition 3.51]{Kollar_Mori_1998}), Theorem~\ref{thm: extending involution} implies that $\mu^0$ extends to an involution $\mu$ on the whole pair $(D_{\Xbar^c},\Diff_{D_{\Xbar^c}}(\Delta_{\Xbar^c}))$.
Let $H=\langle G,\mu\rangle$ be the subgroup of $\Aut(D_{\Xbar^c}^n/T,\Diff_{D_{\Xbar^c}^n}(\Delta_{\Xbar^c}))$ generated by $G$ and $\mu$, which is finite as the automorphism group of a relative stable pair is finite by~\cite{Fujino_RelativeLogPluricanonical}.

Since $(D_{\Xbar^c}^n,\Diff_{D_{\Xbar^c}^n}(\Delta_{\Xbar^c}))$ is the log canonical model of $(D_{\Xbar^m}^n,\Diff_{D_{\Xbar^m}^n}(\Delta_{\Xbar^m}))$ and the morphism $D_{\Xbar^m}^n\rightarrow D_{\Xbar^c}^n$ is crepant by Lemma~\ref{lem: conductor for sklt induction}\eqref{point 4}, we get a canonical inclusion
\[
    H\subset\Bim(D_{\Xbar^m}^n/T,\Diff_{D_{\Xbar^m}^n}(\Delta_{\Xbar^m}))
\]
by Remark~\ref{rmk: Bir is a crepant invariant}.
By the inductive assumption and Lemma~\ref{lem: Artin-Tate}, the inclusions
\begin{equation}\label{eq: finite inclusion conductor}
\begin{tikzcd}
    \widetilde{R}_D^H:=R_{D_{\Xbar^m}^n}^H\cap R_{D_{\Xbar^m}}\arrow[r,hookrightarrow] & R_{D_{\Xbar^m}} \arrow[r,hookrightarrow] & R_{D_{\Xbar^m}^n}
\end{tikzcd}
\end{equation}
are finite extensions. Now, recall that $R_{\Xbar}\simeq R_{\Xbar^m}$, and by the inclusions in Lemma~\ref{lem: conductor for sklt induction}\eqref{point 5} we have a commutative diagram
\begin{equation}\label{diag: restriction}
\begin{tikzcd}
    R_{D_{\Xbar}^n}\arrow[r,hookleftarrow] & R_{D_{\Xbar^m}^n}\arrow[r,phantom,"\simeq"] & R_{D_{\Xbar^c}^n}\\
    R_{\Xbar}\arrow[r,phantom,"\simeq"]\arrow[u,"\res_n"] & R_{\Xbar^m}\arrow[u,"\res_n"]\arrow[r,phantom,"\simeq"] & R_{\Xbar^c}\arrow[u,"\res_n"]
\end{tikzcd}
\end{equation}
showing that
\[
    R_X=
    \res_n^{-1}(R_{D_{\Xbar^m}^n}^{\mu})\supset\res_n^{-1}(R_{D_{\Xbar^m}^n}^H)=\res^{-1}(\widetilde{R}_D^H),
\]
where we have used Lemma~\ref{lem: log canonical ring slc pair} and Remark~\ref{rmk: different invariants mu}. Since the inclusions in~\eqref{eq: finite inclusion conductor} are finite extensions, they are in particular integral. By Lemma~\ref{lem: surjectivity restriction minimal model}, the restriction map between $\ell$-th Veronese subalgebras
\[
\begin{tikzcd}
    \res\colon R_{\Xbar^m}^{(\ell)}\arrow[r] & R_{D_{\Xbar^m}}^{(\ell)}
\end{tikzcd}
\]
is surjective for $\ell$ sufficiently divisible. In particular, the pullback of~\eqref{eq: finite inclusion conductor}

\begin{equation}\label{eq: extension}
\begin{tikzcd}
    \res^{-1}(\widetilde{R}_D^H)\arrow[r,hookrightarrow] & R_{\Xbar^m}
\end{tikzcd}
\end{equation}

induces an integral extension between appropriate Veronese subalgebras, by Lemma~\ref{lem: pullback integral extension}. Since any graded algebra is integral over any of its Veronese subalgebras, the extension in~\eqref{eq: extension} is also an integral extension, hence finite, as $R_{\Xbar^m}$ is finitely generated over $T$.

Since $\res^{-1}(\widetilde{R}_D^H)\subset R_X\subset R_{\Xbar}\simeq R_{\Xbar^m}$, the extension $R_X\subset R_{\Xbar}$ is also finite. In particular, by Lemma~\ref{lem: Artin-Tate}, $(X,\Delta_X)$ has finitely generated log canonical algebra, thus admitting a good stable model by Theorem~\ref{thm: projRX is stable model}. Note that there could be vertical irreducible components of $X$ over $T$, but they have image in $T^0$ by assumption, where we already know $\Projan_T R_X$ to be the stable model again by assumption.

We also have inclusions $R_X\cap R_{\Xbar}^G\subset\res^{-1}(\widetilde{R}_D^H)\subset R_{\Xbar}$. By Lemma~\ref{lem: Artin-Tate}, to prove the theorem it is then enough to show that the inclusion $R_X\cap R_{\Xbar}^G\subset\res^{-1}(\widetilde{R}_D^H)$ is integral, for which we can pass to Veronese subalgebras as needed. Note that $R_X\cap R_{\Xbar}^G\subset\res^{-1}(\widetilde{R}_D^H)$ is not an equality: $\alpha \in R_X^{(\ell)}$ maps to $\widetilde{R}_D^H$ if its image $\res(\alpha)$ in $R_{D_{\Xbar^m}^n}$ is invariant under $G$, but $\alpha$ could still fail to be $G$-invariant as a section of $R_{\Xbar}$, as of course the restriction map $\res_n$ is not injective.

Since $G$ is linearly reductive, the map
\[
\begin{tikzcd}
    \mathrm{res}^G\colon (R_{\Xbar}^{(\ell)})^G\simeq (R_{\Xbar^m}^{(\ell)})^G \arrow[r] & (R_{D_{\Xbar^m}}^{(\ell)})^G
\end{tikzcd}
\]
is again surjective for $\ell$ sufficiently divisible. Let $I$ be the kernel ideal of the restriction morphism $\res$, hence $I^G=I\cap R_{\Xbar^m}^G$ is the kernel of $\res^G$.

Let $\alpha\in\res^{-1}(\widetilde{R}_D^H)$ have even, sufficiently divisible degree, and
\[
    \bar{\alpha} = \mathrm{res}(\alpha) \in R_{D_{\Xbar^m}}\cap R_{D_{\Xbar^m}^n}^H \subset R_{D_{\Xbar^m}}^G
\]
its image. By surjectivity of $\res^G$, there exists $\beta\in R_{\Xbar}^G$ s.t. $\res(\beta) = \bar{\alpha}$.
In particular, $\mathrm{res}(\beta) \in R_{D_{\Xbar^m}}\cap R_{D_{\Xbar^m}^n}^\mu$, hence $\beta \in R_X\cap R_{\Xbar}^G$. By construction, $\alpha-\beta\in I$.
Therefore, to conclude it is enough to show that every $\alpha \in I$ of sufficiently divisible degree is integral over $R_X\cap R_{\Xbar}^G$.

Note that $I \subset R_X$, as the zero-function on $R_{D_{\Xbar}^n}$ is trivially $\mu$-invariant.
Consider the monic polynomial $P(X) := \prod_{g \in G}(X - g \cdot \alpha)$, hence $P(\alpha) = 0$ (look at $g = \text{identity}$).
The polynomial is $G$-invariant by construction, with coefficients being (up to sign) symmetric functions in $g \cdot \alpha$ for $g \in G$ varying. In particular, the coefficients are in $I^G \subset R_X\cap R_{\Xbar}^G$, proving the result.
\end{proof}

\subsection{Proof of Theorem~\ref{thm: general version}}\label{subsec: proof main theorem}

The idea of the proof of Theorem~\ref{thm: general version} is to perturb the semi-dlt pair $(X,\Delta_X)$ by adding $\epsilon A_X$ for a relatively ample divisor $A_X$ over $T$, with $\epsilon>0$ sufficiently small. The resulting pair is equivalent to a semi-klt pair, so that we can apply Theorem~\ref{thm: equivariant sklt} to obtain a stable model $\phi^c_+\colon X\dashrightarrow X^c_+$ that does not depend on $\epsilon$ for $0<\epsilon\ll1$. One difficulty is that there is a jumping behavior at $\epsilon=0$, where the canonical model can be different, if it exists. To deduce the finite generation of $R_X=R(X/T,K_X+\Delta_X)$, we show that $K_{X^c_+}+(\phi^c_+)_*\Delta_X$ is a semi-ample $\QQ$-Cartier divisor over $T$. As usual, we will tacitly restrict $T$ when needed, as explained in Remark~\ref{rmk: reductions}.
%

To implement the perturbation strategy, the first step consists in showing that the perturbed pair $(X,\Delta_X+\epsilon A_X)$ satisfies the same assumptions as Theorem~\ref{thm: equivariant sklt}. We start by showing that we can replace a perturbed pair with a semi-klt one that is well-behaved with respect to any given bimeromorphic map.

\begin{lemma}\label{lem: improving boundary by perturbation}
    Let $(X,\Delta_X)$ be a semi-dlt pair, $f\colon X\rightarrow T$ be a projective morphism to a complex analytic space $T$, and $A_X$ an $f$-ample $\QQ$-Cartier divisor. Let $\phi\colon X\dashrightarrow Y$ be a bimeromorphic map over $T$, whose exceptional locus does not contain any stratum of $D_X$. Then, after restricting to any relatively compact Stein open subset of $T$, there exists an integer $c>0$ and a $\QQ$-divisor $B_X$ such that for all rational $0<\epsilon\ll1$ the following holds:
    \begin{enumerate}
        \item\label{req 1} on the open locus of $X$ where $X$ is SNC, we have that $B_X$ is $\QQ$-Cartier and $f$-ample, and $\Supp B_X$ is SNC-adapted;
        \item\label{req 2} $K_X+\Delta_X+\epsilon(cA_X)\sim_{\QQ,T}K_X+(1-\epsilon)\Delta_X+\epsilon B_X$;
        \item\label{req 3} $(X,(1-\epsilon)\Delta_X+\epsilon B_X)$ is semi-klt, in particular $\lfloor\epsilon B_X\rfloor=0$;
        \item\label{req 4} the exceptional locus of $\phi$ does not contain any component of $\Supp B_X$ and $\Supp B_X\cap V$ for any stratum $V$ of $D_X$ for which $\Supp B_X\cap V\not=\emptyset$.
    \end{enumerate}
\end{lemma}
\begin{proof}
    We follow the proof of~\cite[Proposition 2.43]{Kollar_Mori_1998} and adapt it our setting. In order to simplify the notation and using Remark~\ref{rmk: reductions}, in the following we will assume that $f:X\rightarrow T$ is obtained from another projective morphism by restriction to a relatively compact Stein open subset.
    
    By Remark~\ref{rmk: existence KX} there is a canonical divisor $K_X$. By~\cite[Theorem 2.12]{Das_Hacon_Pau_4dimensionKahler}, we can choose positive integers $m$ and $m'$ such that $m\Delta_X$ is an integral divisor, $m'A_X$ is Cartier and $f$-very-ample, and the sheaf $\cO_X(m\Delta_X+m'A_X)$ is generated by global sections over $T$.
    Set $c=m'/m$, where we can assume that $m$ divides $m'$. Let $B_X'\in|m\Delta_X+m'A_X|_T$ be a general member, and let $B_X=(1/m)B_X'$.
    Then, $-m\Delta_X+mB_X\sim_T m'A_X$ is Cartier, thus
    \[
        K_X+(1-\epsilon)\Delta_X+\epsilon B_X\sim_{\QQ,T} K_X+\Delta_X+\epsilon(c A_X)
    \]
    and is $\QQ$-Cartier. Let $Z\subset X$ be the analytic subset in the definition of semi-dlt pair (Definition~\ref{def: dlt/sdlt}), that is, the pair $(U=X\setminus Z,\Delta_X|_U)$ is semi-SNC, and the discrepancies over any center in $Z$ is strictly greater than $-1$. Note that $U$ is contained in the larger open $U'$ where $X$ is SNC, and in particular Gorenstein. Note that $U'$ is a Zariski-open subset by Remark~\ref{rmk: analycity SNC locus}. On $U'$ the divisor $\Delta_X$ is $\QQ$-Cartier, thus the same holds for $B_X$. Choosing $m'$ big enough, we can assume that $B_X$ is again $T$-ample on $U'$. Moreover, by base-point freeness and~\cite[Theorem 2.20]{Das_Hacon_Pau_4dimensionKahler} applied to the normalization of $U$, we can choose general $B_X$ so that the support $\Supp B_X$ has normal crossing with $D_X\cup\Supp\Delta_X$ on $U$.
    Then, for $0\leq\epsilon\ll1$ we have that the perturbed pair is semi-dlt on $U$, hence semi-klt as soon as $\lfloor\epsilon B_X\rfloor=0$. On $Z$, any sufficiently small perturbation stays slc and does not create additional log canonical centers, thus $(X,(1-\epsilon)\Delta_X+\epsilon B_X)$ is everywhere semi-klt for $0\leq\epsilon\ll1$.

    We are left with showing that $\Ex(\phi)$ does not contain any component of $B_X$ and $B_X\cap V$ for any stratum $V$ of $D_X$. On the other hand, this follows from the global generation of $\cO_X(m\Delta_X+m'A_X)$ over $T$, as we just need to choose $B_X$ so that it contains no irreducible component of $\Ex(\phi)$ and $\Ex(\phi)\cap V$.
\end{proof}

The following lemma shows that the perturbed pair still admits a stable model over $T^0$.

\begin{lemma}\label{lem: existence stable model perturbed pair}
    Let $(X,\Delta_X)$ be a semi-dlt pair and $f\colon X\rightarrow T$ a projective morphism to a complex analytic space $T$. Suppose the pair admits a good stable model $(X^c,\Delta_{X^c})$ over $T$, and that the exceptional locus of $\phi^c\colon X\dashrightarrow X^c$ does not contain any stratum of the conductor $D_X$. Let $A_X$ be an $f$-ample $\QQ$-Cartier divisor on $X$. Then, the pair $(X,\Delta_X+\epsilon A_X)$ admits a good stable model $\phi_{\epsilon}^c\colon X\dashrightarrow X_{\epsilon}^c$ for all $0<\epsilon\ll 1$, and $\Ex(\phi_{\epsilon}^c)$ does not contain any stratum of $D_X$.
\end{lemma}
\begin{proof}
    The statement is local on $T$, therefore we can restrict $T$ to a relatively compact Stein open subset, work around a Stein compact subset $W$ of of $T$ satisfying $\cP$, and further restrict $T$ around $W$ when necessary to apply results from the minimal model program, as in Remark~\ref{rmk: reductions}. Then, by Lemma~\ref{lem: improving boundary by perturbation}, there is an integer $c>0$ and a $\QQ$-divisor $B_X$ such that $K_X+\Delta_X+\epsilon(cA_X)\sim_{\QQ,T}K_X+(1-\epsilon)\Delta_X+\epsilon B_X$ and $(X,(1-\epsilon)\Delta_X+\epsilon B_X)$ is semi-klt for all $0<\epsilon\ll1$. In particular, the normalization pair $(\Xbar,D_{\Xbar}+(1-\epsilon)\Delta_{\Xbar}+\epsilon B_{\Xbar})$ is dlt.
    
    Moreover, the normalization pair admits a log canonical model by the complex analytic version~\cite{Fujino_MMPanalytic} of the main results of~\cite{BCHM}, and~\cite[Proposition 2.43]{Kollar_Mori_1998}.
    Eventually after taking a small $\QQ$-factorial modification over $W$ (\cite[Theorem 1.24]{Fujino_MMPanalytic}), the MMP with scaling of the ample divisor $A_{\Xbar}:=\pi_X^*A_X$ terminates and the resulting $\QQ$-factorial minimal model $\phibar^m\colon \Xbar\dashrightarrow\Xbar^m$ of $(\Xbar/T,D_{\Xbar}+\Delta_{\Xbar}+\epsilon(cA_{\Xbar}))$ around $W$ does not depend on $0\leq\epsilon\ll1$, see Lemma~\ref{lem: existence minimal model if log canonical} and~\cite[\S2]{Hacon_Xu_dslt}. Note in particular that the minimal model is the same as that of $(\Xbar,D_{\Xbar}+\Delta_{\Xbar})$. We set $B_{\Xbar^m}:=\phibar^m_*B_{\Xbar}$. Let us also denote by $R_{D_{\Xbar^m}}^{\epsilon}$ the log canonical algebra of the perturbed pair $(D_{\Xbar^m},\Diff_{D_{\Xbar^m}}(\Delta_{\Xbar^m})+\epsilon B_{\Xbar^m}|_{D_{\Xbar^m}})$, and similarly for all other perturbed pairs.
    
    Since $(X,\Delta_X)$ admits a stable model $\phi^c\colon X\dashrightarrow X^c$ over $T$, there is an involution on the pair $(D_{\Xbar^c}^n,\Diff_{D_{\Xbar^c}^n}(\Delta_{\Xbar^c}))$ that defines the gluing relations. Since $\phibar_{D^n}\colon D_{\Xbar^m}^n\rightarrow D_{\Xbar^c}^n$ is bimeromorphic by the assumption on the exceptional locus of $\phi^c$, and the morphism is crepant (see the proof of Lemma~\ref{lem: birational coeff 1 boundary then stable model} or~\cite[Proposition 4.6]{Kollar-singularities}), $\phibar_{D^n}$ is a crepant log canonical model over $T$. By Remark~\ref{rmk: Bir is a crepant invariant}, the involution on $D_{\Xbar^c}^n$ induces a B-bimeromorphic involution $\mu\in\Bim(D_{\Xbar^m}^n/T,\Diff_{D_{\Xbar^m}^n}(\Delta_{\Xbar^m}))$. Since $(D_{\Xbar^m}^n,\Diff_{D_{\Xbar^m}^n}(\Delta_{\Xbar^m}))$ is a minimal model, $\mu$ is an holomorphic automorphism on an open $U\subset D_{\Xbar^m}^n$ whose complement is of codimension $\geq2$, see for instance~\cite[Theorem 3.52, Corollary 3.54]{Kollar_Mori_1998}. In particular, to check that $\mu\in\Bim(D_{\Xbar^m}^n/T,M)$ for some boundary $M$, it is enough to check so for $\mu|_U$.

    By Lemma~\ref{lem: improving boundary by perturbation}\eqref{point 4} applied to $\phi^c$, none of the components of the divisor $B_{\Xbar}|_{D_{\Xbar}}$ is contained in the exceptional locus of the bimeromorphic map $\phibar^m_D\colon D_{\Xbar}\dashrightarrow D_{\Xbar^m}$. In particular, $\phibar_D^m$ is an isomorphism at the generic points of $B_{\Xbar}|_{D_{\Xbar}}$ and maps it bimeromorphically onto $B_{\Xbar^m}|_{D_{\Xbar^m}}$. Since the diagram
    \[
    \begin{tikzcd}
        D_{\Xbar}^n\arrow[r,dashed]\arrow[d,"\simeq","\text{gluing}"'] & D_{\Xbar^m}^n\arrow[d,dashed,"\mu"]\\
        D_{\Xbar}^n\arrow[r,dashed] & D_{\Xbar^m}^n
    \end{tikzcd}
    \]
    is commutative and $\mu$ is a regular automorphism in codimension 1, it follows that
    \[
        \mu\in\Bim(D_{\Xbar^m}^n/T,\Diff_{D_{\Xbar^m}^n}(\Delta_{\Xbar^m})+\epsilon B_{\Xbar^m}|_{D_{\Xbar^m}}).
    \]
    By Lemma~\ref{lem: conductor for sklt induction}\eqref{point 3}, the pair $(D_{\Xbar^m},\Diff_{D_{\Xbar^m}}(\Delta_{\Xbar^m})+\epsilon B_{\Xbar^m}|_{D_{\Xbar^m}})$ is semi-klt, and it admits a good stable model over $T$ by Lemma~\ref{lem: birational coeff 1 boundary then stable model}. Therefore, we can apply Theorem~\ref{thm: equivariant sklt} to obtain that the inclusions
    \[
    \begin{tikzcd}
        R_{D_{\Xbar^m}}^{\epsilon}\cap (R_{D_{\Xbar^m}^n}^{\epsilon})^{\langle\mu\rangle}\arrow[r,hookrightarrow] & R_{D_{\Xbar^m}}^{\epsilon}\arrow[r,hookrightarrow] & R_{D_{\Xbar^m}^n}^{\epsilon}
    \end{tikzcd}
    \]
    are finite extensions of $\cO_T$-algebras, where we are using Lemma~\ref{lem: Artin-Tate} and the convention of Definition~\ref{def: invariants under Bir} (see also Remark~\ref{rmk: different invariants mu}). By Lemma~\ref{lem: surjectivity restriction minimal model}, the restriction morphism $\res\colon R_{\Xbar^m}^{\epsilon}\rightarrow R_{D_{\Xbar^m}}^{\epsilon}$ induces a surjection between $\ell$-Veronese subalgebras for $\ell$ sufficiently divisible. By Lemma~\ref{lem: pullback integral extension}, it follows that
    \[
    \begin{tikzcd}
        \res^{-1}(R_{D_{\Xbar^m}}^{\epsilon}\cap (R_{D_{\Xbar^m}^n}^{\epsilon})^{\langle\mu\rangle})\arrow[r,hookrightarrow] & R_{\Xbar}^{\epsilon}
    \end{tikzcd}
    \]
    is also an integral extension, hence finite as $R_{\Xbar}^{\epsilon}$ is finitely generated. By Lemma~\ref{lem: log canonical ring slc pair} and Remark~\ref{rmk: different invariants mu}, the inverse image on the left is identified with a Veronese subalgebra of $R_X^{\epsilon}$ of even degree. Therefore, $R_X^{\epsilon}\hookrightarrow R_{\Xbar}^{\epsilon}$ is also a finite extension, hence $R_X^{\epsilon}$ is finitely generated for all $0\leq\epsilon\ll1$ by Lemma~\ref{lem: Artin-Tate}.

    Since $\phi^c$ is an isomorphism at the generic point of any stratum of $D_X$, and $\phi_{\epsilon}^c$ is induced by twisting with an $f$-ample Cartier divisor, the same property holds for $\phi_{\epsilon}^c$. In particular, it is a good stable model by Theorem~\ref{thm: projRX is stable model}.
\end{proof}

We also need to choose the perturbation so that it does not create new log centers on $D_X$ that map to $T\setminus T^0$. For this, we will use the following simple lemma characterizing centers of SNC pairs whose minimal log discrepancy is 1 and that do not dominate $T$.

\begin{lemma}\label{lem: mld=1}
    Let $Y$ be a connected complex manifold and let $\Delta_Y=\sum_{i\in I}a_i\Delta_i$ be an SNC divisor such that $a_i\leq 1$ and are non-zero for all $i\in I$. Let $\Delta_Y^{>0}=\sum_{i\in I:a_i>0}a_i\Delta_i$ be the positive part of $\Delta_Y$, and $f\colon Y\rightarrow T$ a projective surjective morphism to a complex analytic space, such that all strata of $\Delta_Y^{>0}$ surject onto $T$.
    Let $V\subset Y$ be an irreducible reduced analytic subspace, whose image $f(V)=Z\subsetneq T$ is a proper subspace of $T$. If $\mathrm{mld}(V,Y,\Delta_Y)=1$, then one of the following is satisfied:
    \begin{enumerate}
       \item\label{discrepancy 0: case 1} There exists $J\subset I$ such that $V$ is an irreducible component of $\cap_{j\in J}\Delta_j$, and $\sum_{j\in J}a_j=|J|-1$. Moreover, in this case there exists $j\in J$ for which $a_j<0$.
       \item\label{discrepancy 0: case 2} There exists $J\subset I$ with $a_j=1$ for all $j\in J$, such that $V$ is an irreducible component of $(\cap_{j\in J}\Delta_j)\cap f^{-1}(Z)$ and has codimension 1 in $(\cap_{j\in J}\Delta_{j})$.
    \end{enumerate}
    In particular, there exists at most a finite number of such centers $V$.
\end{lemma}
\begin{proof}
    Let $J$ be the subset of indices $j\in I$ for which $V\subset\Delta_j$. By~\cite[2.10]{Kollar-singularities}, we know that
    \begin{equation}\label{eq: equality mld}
        1=\mathrm{mld}(V,Y,\Delta_Y)=\codim_{Y}V-\sum_{j\in J}a_j\geq\codim_{Y}V-|J|,
    \end{equation}
    as $a_j\leq1$. In particular, $|J|+1\geq\codim_{Y}V\geq|J|$.
    
    If $\codim_{Y}V=|J|$, then $V$ is an irreducible component of $\cap_{j\in J}\Delta_j$, hence a stratum of $\Delta_Y$. Since $f(V)\subsetneq T$ and all strata of $\Delta_Y^{>0}$ dominate $T$, it follows that there exists $j\in J$ for which $a_j<0$. The equality $\sum_{j\in J}a_j=|J|-1$ follows from~\eqref{eq: equality mld}.
    
    Suppose now that $\codim_{Y}V=|J|+1$. From~\eqref{eq: equality mld} and the assumption $a_j\leq1$, we get that  $a_j=1$ for all $j\in J$, hence $\Delta_j$ is a component of $\Delta_Y^{>0}$ for all $j\in J$. By assumption, $V\subseteq(\cap_{j\in J}\Delta_j)\cap f^{-1}(Z)$, and $\cap_{j\in J}\Delta_j$ dominates $T$. Therefore, $V$ is an irreducible component of $(\cap_{j\in J}\Delta_j)\cap f^{-1}(Z)$ and has codimension 1 in $\cap_{j\in J}\Delta_j$.
\end{proof}

%

We are ready to prove the main result of the paper.

\begin{proof}[Proof of Theorem~\ref{thm: general version}]
    By Theorem~\ref{thm: equivariant sklt}, we know that the statement holds for all semi-klt pairs. By Remark~\ref{rmk: reductions}, we can assume to have restricted $T$ to a relatively compact Stein open subset. We can also assume $X$ to be non-normal, hence $D_X\not=\emptyset$.

    Fix an $f$-ample line bundle $A_X$ on $X$, and denote by $A_{\Xbar}$ its pullback to $\Xbar$. By Lemma~\ref{lem: improving boundary by perturbation} applied to the bimeromorphic map $\phi^c\colon X\dashrightarrow X^{0,c}$, there exists a boundary $\QQ$-divisor $B_X$ such that
    \[
        K_X+(1-\epsilon)\Delta_X+\epsilon B_X\sim_{\QQ,T}K_X+\Delta_X+\epsilon(cA_X)
    \]
    for every rational $0<\epsilon\ll1$, and $(X,(1-\epsilon)\Delta_X+\epsilon B_X)$ is semi-klt. Lemma~\ref{lem: existence stable model perturbed pair} says that the restriction to $X^0$ of the perturbed pair still admits a good stable model $\phi_{\epsilon}^{0,c}\colon X^0\dashrightarrow X_{\epsilon}^{0,c}$ such that $\Ex(\phi_{\epsilon}^{0,c})$ does not contain any stratum of $D_{X^0}$.
    To apply Theorem~\ref{thm: equivariant sklt}, we need to show that we can choose $B_X$ so that $(X,(1-\epsilon)\Delta_X+\epsilon B_X)$ does not admit new log centers on $D_X$ that do not intersect $X^0$. Since the discrepancies are rational with denominator bounded by the index (\cite[Definition 2.49]{Kollar-singularities}) of the log canonical divisor and $0<\epsilon\ll1$, it is enough to show that we can choose $B_X$ so that its support does not contain centers $V$ lying on $D_X$ with $\mathrm{mld}(V,X,\Delta_X)=1$ and mapping to $T\setminus T^0$. By taking a log resolution of the normalization of $(X,\Delta_X)$ and applying Lemma~\ref{lem: mld=1} to a neighborhood of $D_{\Xbar}$ in each connected component of $\Xbar$, we see that such centers $V$ are in finite number. In particular, in Lemma~\ref{lem: improving boundary by perturbation} we can choose $B_X$ not containing them.
    
    Therefore, for all $0<\epsilon\ll1$, we can apply Theorem~\ref{thm: equivariant sklt} to produce good stable models
    \[
    \begin{tikzcd}
        \phi_{\epsilon}^c\colon (X,(1-\epsilon)\Delta_X+\epsilon B_X)\arrow[r,dashed] & (X_{\epsilon}^c,(1-\epsilon)\Delta_{X_{\epsilon}^c}+\epsilon B_{X_{\epsilon}^c}),
    \end{tikzcd}
    \]
    whose normalization is a log canonical model
    \[
    \begin{tikzcd}
        \phibar_{\epsilon}^c\colon (\Xbar,D_{\Xbar}+(1-\epsilon)\Delta_X+\epsilon B_{\Xbar})\arrow[r,dashed] & (\Xbar_{\epsilon}^c,D_{\Xbar_{\epsilon}^c}+(1-\epsilon)\Delta_{\Xbar_{\epsilon}^c}+\epsilon B_{\Xbar_{\epsilon}^c}).
    \end{tikzcd}
    \]
    Moreover, by assumption the dlt pair $(\Xbar,D_{\Xbar}+\Delta_{\Xbar})$ admits a good canonical model $\phibar^c_0=\phibar^c\colon \Xbar\dashrightarrow\Xbar_0^c=\Xbar^c$. By Lemma~\ref{lem: existence minimal model if log canonical} and~\cite[\S2]{Hacon_Xu_dslt}, eventually after taking a small $\QQ$-factorial modification (\cite[Theorem 1.24]{Fujino_MMPanalytic}), the $(K_{\Xbar}+D_{\Xbar}+\Delta_{\Xbar})$-MMP over $T$ with $A_{\Xbar}$-scalings terminates with a bimeromorphic contraction over $T$
    \[
    \begin{tikzcd}
        \phibar^m\colon \Xbar\arrow[r,dashed] & \Xbar^m
    \end{tikzcd}
    \]
    that is a $\QQ$-factorial dlt $K_{\Xbar}+D_{\Xbar}+\Delta_{\Xbar}+\epsilon(cA_{\Xbar})$-good minimal model for all $0\leq \epsilon\ll1$.
    We emphasize that this includes the case $\epsilon=0$. As usual, set $\Delta_{\Xbar^m}=\phibar^m_*\Delta_{\Xbar}$, and similarly for all the other divisors.
    Since $\Xbar^m$ is a good minimal model, for all $0\leq\epsilon\ll1$ there exist morphisms
    \[
    \begin{tikzcd}
        \phibar_{\epsilon}\colon \Xbar^m\arrow[r] & \Xbar_{\epsilon}^c
    \end{tikzcd}
    \]
    to the respective stable models through which $\phibar_{\epsilon}^c$ factors.
    
    We claim that $\phibar_{\epsilon}$ and $\Xbar_{\epsilon}^c$ do not depend on $\epsilon$ for $0<\epsilon\ll1$; note that these could (and usually do) change at $\epsilon=0$.
    To prove this, recall that, by definition, $\phibar_{\epsilon}$ can be characterized as the bimeromorphic contraction that contracts the curves $C$ that are vertical over $T$ and satisfy
    \[
        C\cdot(K_{\Xbar^m}+D_{\Xbar^m}+(1-\epsilon)\Delta_{\Xbar^m}+\epsilon B_{\Xbar^m})=C\cdot(K_{\Xbar^m}+D_{\Xbar^m}+\Delta_{\Xbar^m}+\epsilon(cA_{\Xbar^m}))=0.
    \]
    On the other hand, $C\cdot(K_{\Xbar^m}+D_{\Xbar^m}+\Delta_{\Xbar^m})\geq0$, thus the above equation would imply that either $C\cdot(K_{\Xbar^m}+D_{\Xbar^m}+\Delta_{\Xbar^m})=0$ or $\epsilon(cA_{\Xbar^m})\cdot C<0$. In the second case, choosing $0<\epsilon<\epsilon'\ll1$, we would then get
    \[
       C\cdot(K_{\Xbar^m}+D_{\Xbar^m}+\Delta_{\Xbar^m}+\epsilon'(cA_{\Xbar^m}))<C\cdot(K_{\Xbar^m}+D_{\Xbar^m}+\Delta_{\Xbar^m}+\epsilon(cA_{\Xbar^m}))=0,
    \]
    contradicting the nefness of $K_{\Xbar^m}+D_{\Xbar^m}+\Delta_{\Xbar^m}+\epsilon'(cA_{\Xbar^m})$. The independence of $\phibar_{\epsilon}^c$ from $0<\epsilon\ll1$ follows.
    We denote by
    \[
    \begin{tikzcd}
        \phibar_+\colon \Xbar^m\arrow[r] & \Xbar_+^c
    \end{tikzcd}
    \]
    the canonical morphism equal to $\phibar_{\epsilon}$ for all $0<\epsilon\ll1$.
    
    Of course, in general $\Xbar_0^c\not=\Xbar_+^c$, as one could drop into a smaller face of the nef cone. On the other hand, if $C$ is contracted by $\phibar_+$, then
    \begin{align*}
        0&=C\cdot\left((K_{\Xbar^m}+D_{\Xbar^m}+\Delta_{\Xbar^m}+\epsilon(cA_{\Xbar^m}))-(K_{\Xbar^m}+D_{\Xbar^m}+\Delta_{\Xbar^m}+\epsilon'(cA_{\Xbar^m})\right))\\
        &=(\epsilon-\epsilon')C\cdot A_{\Xbar^m},
    \end{align*}
    for all $0<\epsilon<\epsilon'\ll1$. Therefore, $C\cdot(K_{\Xbar^m}+D_{\Xbar^m}+\Delta_{\Xbar^m})=0$, and $C$ is contracted by $\phibar_0$ as well. In particular, $\phibar_0$ factors through a bimeromorphic morphism
    \[
    \begin{tikzcd}
        \psi\colon \Xbar_+^c\arrow[r] & \Xbar_0^c,
    \end{tikzcd}
    \]
    see for instance~\cite[Lemma 1.15]{Debarre-higher_dimensional_AG}.
    Note that the log canonical model $\phibar_+^c:=\phibar_+\circ\phibar^m\colon \Xbar\dashrightarrow\Xbar_+^c$ does not contract any stratum of the conductor $D_{\Xbar}$. Indeed, otherwise $\phibar^c=\psi\circ\phibar_+^c$ would also contract some stratum of $D_{\Xbar}$, contradicting the assumption that $\phi^{0,c}$ does not contract any stratum of the conductor $D_X$ of $X$.
    
    We claim that $K_{\Xbar_+^c}+D_{\Xbar_+^c}+\Delta_{\Xbar_+^c}$ is semi-ample over $T$.
    To prove this, it is enough to show that $\psi^*(K_{\Xbar_0^c}+D_{\Xbar_0^c}+\Delta_{\Xbar_0^c})\sim_{\QQ,T}K_{\Xbar_+^c}+D_{\Xbar_+^c}+\Delta_{\Xbar_+^c}$. On the other hand,
    \begin{align*}
        \psi^*(K_{\Xbar_0^c}+D_{\Xbar_0^c}+\Delta_{\Xbar_0^c}) & \sim_{\QQ,T} \phi_{+*}\phi_+^*\psi^*(K_{\Xbar_0^c}+D_{\Xbar_0^c}+\Delta_{\Xbar_0^c})\\
        & =\phi_{+*}\phi_0^*(K_{\Xbar_0^c}+D_{\Xbar_0^c}+\Delta_{\Xbar_0^c})\\
        & = \phi_{+*}(K_{\Xbar^m}+D_{\Xbar^m}+\Delta_{\Xbar^m})\\
        & \sim_{\QQ,T} K_{\Xbar_+^c}+D_{\Xbar_+^c}+\Delta_{\Xbar_+^c}
    \end{align*}
    where the first $\QQ$-linear equivalence follows from the projection formula, and the last by definition. This proves that $K_{\Xbar_+^c}+D_{\Xbar_+^c}+\Delta_{\Xbar_+^c}$ is semi-ample over $T$.

    Now, we work at the demi-normal level, and we claim that also $\phi_{\epsilon}^c\colon X\dashrightarrow X_{\epsilon}^c$ is independent of $\epsilon$ as long as $0<\epsilon\ll1$. From what we have proved above, the normalization $\phibar_{+}^c\colon \Xbar\dashrightarrow\Xbar_+^c$ of $\phi_{\epsilon}^c$ is independent of $\epsilon$ in the same range, and it restricts to a bimeromorphic map between conductors, thus also the meromorphic involution on $D_{\Xbar_+^c}$ is independent of $\epsilon$. Since a demi-normal pair is uniquely determined by its normalization triple (see \S\ref{subsec: normal and demi-normal modification} or \cite[Proposition 5.3]{Kollar-singularities}), it follows that $X\dashrightarrow X_{\epsilon}^c$ is constant for $0<\epsilon\ll1$, as wanted. As before, we denote $X_{\epsilon}^c$ by $X_+^c$ for $0<\epsilon\ll1$, and similarly for the relevant divisors.
    
    Recall that we also know
    \[
    \begin{tikzcd}
        \pi_+\colon (\Xbar_+^c,D_{\Xbar_+^c}+(1-\epsilon)\Delta_{\Xbar_+^c}+\epsilon B_{\Xbar_+^c})\arrow[r] & (X_+^c,(1-\epsilon)\Delta_{X_+^c}+\epsilon B_{X_+^c})
    \end{tikzcd}
    \]
    to be the normalization. It follows that
    \[
        \pi_+^*(K_{X_+^c}+\Delta_{X_+^c})\sim_{\QQ,T}K_{\Xbar_+^c}+D_{\Xbar_+^c}+\Delta_{\Xbar_+^c},
    \]
    which is semi-ample. By~\cite[Theorem 1.1]{Fujino_RelativeLogPluricanonical}, we obtain that $K_{X_+^c}+\Delta_{X_+^c}$ is semi-ample as well. Therefore, the associated $\cO_T$-algebra $R(X_+^c/T,K_{X_+^c}+\Delta_{X_+^c})\simeq R_X$ is locally finitely generated~\cite[Lemma 2.36]{Fujino_MMPanalytic}. Note that~\cite[Lemma 2.36]{Fujino_MMPanalytic} assumes $X$ to be normal, but the proof works in the same way replacing~\cite[Lemma 2.26]{Fujino_MMPanalytic} with Lemma~\ref{lem: veronese finite generation}. By Theorem~\ref{thm: projRX is stable model}, $X^c:=\Projan_T R_X$ with the induced boundary is the good stable model of $(X,\Delta_X)$. Note that there could be vertical irreducible components of $X$ over $T$, but they have image in $T^0$, where we already know $\Projan_T R_X$ to be the good stable model by assumption.
\end{proof}

\subsection{Proof of corollaries}\label{subsec: corollaries}

In this subsection, we prove Corollaries~\ref{cor: Hacon-Xu generalization} and~\ref{cor: analog Birkar} in this order.

\begin{proof}[Proof of Corollary~\ref{cor: Hacon-Xu generalization}]
    To prove Corollary~\ref{cor: Hacon-Xu generalization}, it is enough to show that $(\Xbar,D_{\Xbar}+\Delta_{\Xbar})$ admits a log canonical model over $T$, and then apply Theorem~\ref{thm: general version}. Since the result is local on $T$, we can shrink $T$ as needed, as explained in Remark~\ref{rmk: reductions}. Then, the result follows from Theorem~\ref{thm: Hacon-Xu analytic}.
\end{proof}

\begin{proof}[Proof of Corollary~\ref{cor: analog Birkar}]
    As explained in Remark~\ref{rmk: reductions}, in the proof we will tacitly restrict $T$ as needed, and assume that $X$ is not normal, so that $D_X\not=\emptyset$. By~\cite[Theorem 1.2]{Enokizono-Hashizume-analytic-MMP} (see~\cite[Theorem 1.1(3)]{Birkar_existenceflips} or~\cite[Theorem 1.6]{Hacon_Xu_dslt} for the algebraic versions), the normalization pair $(\Xbar,D_{\Xbar}+\Delta_{\Xbar})$ admits both a good minimal model and a log canonical model over $T$. Taking the Stein factorization of $f$, we may further assume that $f$ is bimeromorphic, does not contract any log center lying on $D_X$, and satisfies $f_*\cO_X\simeq\cO_T$. Set $X^0=X\setminus\Ex(f)$ and $T^0=T\setminus f(\Ex(f))$, thus $X^0$ intersects all log centers lying on $D_X$. By connectedness of the fibers, we have $X^0=f^{-1}(T^0)\simeq T^0$.
    Therefore, all assumptions of Theorem~\ref{thm: general version} are satisfied, hence $(X,\Delta_X)$ admits a stable model over $T$.
\end{proof}

\section{Semi-log canonical MMP}\label{sec: slc MMP}

In this section we construct certain MMP steps for projective morphisms between complex analytic spaces with source an slc pair. General MMP steps can fail to exist even in the algebraic case by~\cite[Examples 4, 5]{Ambro-Kollar-slc-MMP} and~\cite[Example 5.4]{fujino_theorems_slc}. As it is probably evident at this point of the manuscript, a sufficient condition for their existence involves log centers lying on the conductor $D_X$. For the algebraic versions of these results, see~\cite{Ambro-Kollar-slc-MMP} and~\cite{Hashizume_sdlt_models}. Our results and theirs appear to be close to optimal in view of the examples mentioned above.

\begin{remark}[Cone and contraction theorem for slc pairs]
    We will use the Fujino's cone and contraction theorem for analytic \emph{quasi-log canonical pairs} \cite[Theorems 9.1, 9.2]{Fujino_QuasiLogAnalytic}. Note that these apply to semi-log canonical pairs $f\colon(X,\Delta_X)\rightarrow T$ by~\cite[Theorem 10.1]{Fujino_QuasiLogAnalytic} and the subsequent discussion. These results are direct generalizations of~\cite{fujino_theorems_slc} to the analytic setting, and provide the usual package that enables us to contract $(K_X+\Delta_X)$-negative extremal rays around Stein compact subsets $W\subset T$ satisfying $\cP$. We remark that the existence of the extremal contraction is only guaranteed after suitably shrinking $T$ around $W$.
\end{remark}

\subsection{MMP for semi-dlt pairs}\label{subsec: sdlt MMP}

We start with establishing the existence of MMP steps for semi-dlt pairs, under some technical conditions and control on the log centers lying on the conductor.

\begin{corollary}\label{cor: existence semidlt MMP}
    Let $(X,\Delta_X)$ be a semi-dlt pair, $f\colon X\rightarrow T$ a projective morphism to a semi-normal Stein space $T$, and $W\subset T$ a Stein compact subset satisfying property $\cP$. Let $F\subset\overline{NE}(X/T;W)$ be a $(K_X+\Delta_X)$-negative extremal ray, and $\varphi=\varphi_F\colon X\rightarrow Z$ the corresponding contraction, defined after suitably shrinking $T$ around $W$. Suppose that $\varphi$ is bimeromorphic and $\Ex(\varphi)\cap D_X$ contains no log centers. Then, $Z$ is demi-normal and there exists an MMP step around $W$
    \[
    \begin{tikzcd}
        (X,\Delta_X)\arrow[rr,dashed,"\psi"]\arrow[rd,"\varphi"'] & & (Y,\Delta_Y)\arrow[ld,"\varphi'"]\\
        & Z
    \end{tikzcd}
    \]
    where $(Y,\Delta_Y)$ is semi-dlt.
\end{corollary}
\begin{proof}
    In the proof, we will tacitly restrict $T$ to a Stein neighborhood of $W$ as needed. Note that $\varphi$ does not contract any component of the conductor. Since $\varphi_*\cO_X\simeq\cO_Z$ by~\cite[Theorem 9.2(4)]{Fujino_QuasiLogAnalytic}, we have that $Z$ is demi-normal by Example~\ref{exm: deminormal modifications}\eqref{exm: deminormal stein factorization}. By~\cite[Theorems 9.1, 9.2]{Fujino_QuasiLogAnalytic}, $Z$ is projective over $T$ and $-(K_X+\Delta_X)\sim_{\QQ,Z}A$ for some $\QQ$-Cartier $\QQ$-divisor $A$ that is $\varphi$-ample over $W$ by Kleiman's ampleness criterion~\cite[Theorem 11.5]{Fujino_ConeContractionAnalytic}. By ampleness and~\cite[Theorems 2.12, 2.20]{Das_Hacon_Pau_4dimensionKahler}, we can choose $A$ so that $(X,\Delta_X+A)$ is slc. By Corollary~\ref{cor: analog Birkar}, there exists a good stable model $\psi\colon (X,\Delta_X)\dashrightarrow(Y,\Delta_Y)$ over $Z$ satisfying Definition~\ref{def:slc-MMP steps}\eqref{slc_step:1}-\eqref{slc_step:2}, and whose normalization $\psibar$ is the log canonical model of $(\Xbar,D_{\Xbar}+\Delta_{\Xbar})$ over $Z$. Moreover, we know from the classical theory that $\psibar$ is an MMP step (\cite{Fujino_MMPanalytic,Enokizono-Hashizume-analytic-MMP}), hence the same holds for $\psi$ by Remark~\ref{rmk: normalization MMP step}. Since $(Y,\Delta_Y)$ is dlt away from the conductor and being dlt is preserved by any MMP step, to prove that $(Y,\Delta_Y)$ is semi-dlt it is enough to show that $\psi$ is an isomorphism at the generic point of every log canonical center of $(X,\Delta_X)$ lying on $D_X$. This follows from the stronger assumption that $\Ex(\varphi)\cap D_X$ contains no log centers.
\end{proof}

\begin{remark}\label{rmk: MMP with scaling}
    If $H$ is a general $f$-ample $\QQ$-Cartier divisor on $X$, $\lambda$ the associated (rational) nef threshold, and $F$ is $K_X+\Delta_X+\lambda H$-trivial, then $K_Y+\Delta_Y+\lambda H_Y$ is again nef and $\QQ$-Cartier, where $H_Y=\psi_*H$; see~\cite[Remark 3.2]{Hashizume_sdlt_models} and~\cite[Theorem 9.2(4.iii)]{Fujino_QuasiLogAnalytic}. In particular, we can run an MMP with scalings of $H$ as long as the conditions on log centers are preserved after each step. For instance, this happens whenever the relative stable base locus of $K_X+\Delta_X$ over $T$ does not contain log centers, as $\Ex(\varphi)$ is contained in it. Moreover, if $K_X$ is $\QQ$-Cartier then the same holds for $K_Y$, by the same argument as in~\cite[Theorem 3.9]{Enokizono-Hashizume-analytic-MMP}.
\end{remark}

\subsection{Semi-dlt modifications}\label{subsec: sdlt modification}

In this subsection we prove the existence of crepant semi-dlt modifications in the analytic setting.
First, we need a semi-resolution lemma.

\begin{lemma}[{\cite[Lemma 4.1]{Hashizume_sdlt_models}}]\label{lem: semiresolution}
    Let $(X,\Delta_X)$ be an slc pair such that every irreducible component of $X$ is normal in codimension one, and $f\colon X\rightarrow T$ a projective morphism to a complex analytic space. Then, after restricting to a relatively compact Stein open subset of $T$, there is a projective bimeromorphic morphism $\Phi\colon\widetilde{X}\rightarrow X$ over $T$, with $\widetilde{X}$ projective over $T$, such that if we define a $\QQ$-divisor $\Xi$ on $\widetilde{X}$ by $K_{\widetilde{X}}+\Xi=\Phi^*(K_X+\Delta_X)$, then the following properties hold.
    \begin{enumerate}
        \item\label{resolution point: 1} $(\widetilde{X},\Xi)$ is an analytic globally embedded simple normal crossing pair. Equivalently, there exists a complex manifold $M$ of dimension $\dim\ \! X+1$ and an SNC $\QQ$-divisor $B$ on $M$ such that $\widetilde{X}+\Supp B$ is an SNC divisor on $M$, $\widetilde{X}$ and $B$ share no common component, and $B|_{\widetilde{X}}=\Xi$.
        \item $\Phi$ is an isomorphism over an open subset of $X$ intersecting every component of the conductor $D_X$.
        \item If we write $\Xi=\Delta_{\widetilde{X}}-E_{\widetilde{X}}$ with $\Delta_{\widetilde{X}}$ and $E_{\widetilde{X}}$ having no common components, then $\Delta_{\widetilde{X}}$ and $E_{\widetilde{X}}$ are both $\QQ$-Cartier.
        \item every point $x\in\widetilde{X}$ is contained in at most two irreducible components of $\widetilde{X}$.
        \item The restriction of $\Phi$ to the conductor $D_{\widetilde{X}}$ of $\widetilde{X}$ defines a bimeromorphic morphism $\Phi|_{D_{\widetilde{X}}}\colon D_{\widetilde{X}}\rightarrow D_X$.
    \end{enumerate}
\end{lemma}
\begin{proof}
    The proof is the same as in~\cite[Lemma 4.1]{Hashizume_sdlt_models}, using~\cite[Theorem 10.1]{Fujino_QuasiLogAnalytic} instead of~\cite[Theorem 1.2]{fujino_theorems_slc}. Moreover, note that the semi-resolution result~\cite[Theorem 1.8]{B-M-except-I} of Bierstone and Milman works in the analytic setting by working over relatively compact open subsets, see~\cite[Remark 1.6]{B-M-except-I}.
\end{proof}

\begin{proof}[Proof of Theorem~\ref{thm: semidlt modifications}]
    In the proof, we tacitly restrict $T$ around $W$. The argument goes along the same lines as in~\cite[Theorem 1.2]{Hashizume_sdlt_models}, but with some additional technical difficulties. Let $\Phi\colon\widetilde{X}\rightarrow X$ be the projective bimeromorphic morphism as in Lemma~\ref{lem: semiresolution}, and write
    \[
        K_{\widetilde{X}}+\Delta_{\widetilde{X}}=\Phi^*(K_X+\Delta_X)+E_{\widetilde{X}},
    \]
    with $\Delta_{\widetilde{X}}$ and $E_{\widetilde{X}}$ both effective, $\QQ$-Cartier, and without any common components. Note that $(\widetilde{X},\Delta_{\widetilde{X}}+E_{\widetilde{X}})$ is a semi-SNC pair. In particular, any divisor over $\widetilde{X}$ that has center contained in $E_{\widetilde{X}}$ has non-negative discrepancy in $(\widetilde{X},\Delta_{\widetilde{X}})$, hence $\Supp E_{\widetilde{X}}$ contains no log centers. Since $K_{\widetilde{X}}+\Delta_{\widetilde{X}}\sim_{\QQ,X}E_{\widetilde{X}}$, the relative stable base locus over $X$ of the log canonical divisor on $\widetilde{X}$ contains no log centers, and this property is preserved under any step of the MMP. By definition, the exceptional locus of any $(K_{\widetilde{X}}+\Delta_{\widetilde{X}})$-negative contraction of an extremal ray over $X$ is contained in the relative stable base locus, hence it does not contain any log center. We would now like to run an MMP over $X$ around the compact subset $f^{-1}(W)$. A technical difficulty, is that $f^{-1}(W)$ is not Stein compact, hence to apply results of the minimal model program for projective morphisms between complex analytic spaces we are forced to run the MMP over $T$ around $W$. We encountered the same problem in the proof of Lemma~\ref{lem: termination very exceptional MMP}, and we will adopt the same solution. Let $A$ be a general $f$-ample divisor on $X$ such that $A\cdot C>2\ \!\dim\ \! X$ for every curve $C$ that is contracted by $f$; this exists after suitably restricting $T$ around $W$. Let $H$ be a general $\Phi$-ample $\QQ$-Cartier divisor on $\widetilde{X}$ such that $K_{\widetilde{X}}+\Delta_{\widetilde{X}}+\Phi^*A+H$ is ample. By~\cite[Lemma 9.4]{Fujino_MMPanalytic} and Remark~\ref{rmk: normalization MMP step}, a $(K_{\widetilde{X}}+\Delta_{\widetilde{X}}+\Phi^*A)$-MMP step over $T$ with scalings of $H$ is also a $(K_{\widetilde{X}}+\Delta_{\widetilde{X}})$-MMP step over $X$ with scalings of $\Phi^*A+H\sim_{\QQ,X}H$. In particular, every contraction associated with a sequence of such MMP steps contracts no log centers. By Corollary~\ref{cor: existence semidlt MMP} and Remark~\ref{rmk: MMP with scaling}, after restricting $T$ around $W$ at each step, we can construct a sequence of $(K_{\widetilde{X}}+\Delta_{\widetilde{X}}+\Phi^*A)$-MMP steps over $T$ around $W$ with scalings of $H$, where all pairs are semi-dlt at each step. By the discussion above, we can see this sequence as an MMP over $X$ with scalings of $H$. Recall that $K_{\widetilde{X}}+\Delta_{\widetilde{X}}$ is $\QQ$-linearly equivalent over $X$ to the exceptional divisor $E_{\widetilde{X}}$, and this holds at each step of the MMP. Passing to the normalization, Lemma~\ref{lem: termination very exceptional MMP} and its proof show that the MMP terminates with a semi-dlt good minimal model $g\colon (Y,\Delta_Y)\rightarrow X$ of $(\widetilde{X},\Delta_{\widetilde{X}})$ over $X$, and contracts all components in $E_{\widetilde{X}}$.
    In particular, $K_Y+\Delta_Y\sim_{\QQ}g^*(K_X+\Delta_X)$.
    
    By Remark~\ref{rmk: MMP with scaling}, $K_Y$ is $\QQ$-Cartier, and the same argument as at the end of the proof of~\cite[Theorem 1.2]{Hashizume_sdlt_models} shows that $(Y,0)$ is semi-terminal.
\end{proof}

\subsection{MMP for slc pairs}\label{subsec: slc MMP}

We are ready to prove Corollary~\ref{cor: existence slc MMP} in full generality.

\begin{proof}[Proof of Corollary~\ref{cor: existence slc MMP}]
    Since $\varphi$ does not contract any component of the conductor and $\varphi_*\cO_X\simeq\cO_Z$ by~\cite[Theorem 9.2(4)]{Fujino_QuasiLogAnalytic}, we have that $Z$ is demi-normal by \S\ref{subsec: normal and demi-normal modification} and Example~\ref{exm: deminormal modifications}\eqref{exm: deminormal stein factorization}. By uniqueness of stable models together with~\cite[2.41 and Corollary 2.43(2)]{Kollar-singularities} (which work in the analytic setting in the same way), it is enough to prove the statement after taking the Galois double cover of $X$ as in Remark~\ref{rmk: double cover} (see also~\cite[5.23]{Kollar-singularities}). In particular, we can assume the irreducible components of $X$ to be normal in codimension 1. Using again that the result is local on $T$, in the following we can suitably restrict $T$ as needed (see also Remark~\ref{rmk: reductions}). In particular, we can apply Theorem~\ref{thm: semidlt modifications} to obtain a crepant semi-dlt modification $\gamma\colon (X',\Delta_{X'})\rightarrow(X,\Delta_X)$ such that $\gamma$ does not contract any stratum of $D_{X'}$. Let $U=X\setminus\Ex(\varphi)$, $U'=\gamma^{-1}(U)$, and $V=Z\setminus\varphi(\Ex(\varphi))$. By connectedness of the fibers, we have $U=\varphi^{-1}(V)$. In particular, $(X',\Delta_{X'})|_V=(U',\Delta_{U'})$ admits a stable model over $V$, which is isomorphic to $V$ itself, and $U'\rightarrow V$ does not contract any stratum of the conductor of $U'$. Moreover, $R_{X'}\simeq R_X$, hence it is enough to show that $R_{X'}$ is finitely generated. Note that $R_{\Xbar'}\simeq R_{\Xbar}$ is a locally finitely generated $\cO_Z$-algebra, and the induced bimeromorphic map $\Xbar\dashrightarrow\Proj R_{\Xbar}$ is the MMP step associated to the contraction $\Xbar\rightarrow\Zbar$. Moreover, $(X'\setminus U')\cap D_{X'}$ contains no log centers. Indeed, since $\gamma$ is crepant, the image along $\gamma$ of a log center lying in $(X'\setminus U')\cap D_{X'}$ would be contained in $(X\setminus U)\cap D_X=\Ex(\varphi)\cap D_X$, which would contradict the assumptions. Therefore, $R_{X'}\simeq R_X$ is locally finitely generated by Theorem~\ref{thm: general version}, and the stable model $X\dashrightarrow\Proj R_X$ over $Z$ is an MMP step.
\end{proof}

\begin{remark}\label{rmk: slc divisorial contraction}
    With the assumptions of Corollary~\ref{cor: existence slc MMP}, suppose further that the restriction of the normalization $\overline{\varphi}\colon \Xbar\rightarrow\Zbar$ to every connected component of $\Xbar$ is a divisorial contraction. Then, $\overline{\varphi}'\colon \Ybar\rightarrow\Zbar$ is an isomorphism. Since $\varphi$ is an isomorphism at the general points of the conductor, $\varphi'\colon Y\rightarrow Z$ is an isomorphism on an open whose complement has codimension at least 2. Therefore, $\varphi'$ is an isomorphism by \S\ref{subsec: normal and demi-normal modification}, hence $(Z,\varphi_*\Delta_X)$ is slc and $\varphi$ the MMP step.
\end{remark}

\bibliographystyle{amsalpha}
\bibliography{library}

\end{document}